\documentclass[11pt,a4paper,twoside]{report}
\usepackage{graphicx}
\usepackage[T1]{fontenc}
\usepackage[latin1]{inputenc}
\usepackage{amssymb}
\usepackage{latexsym}
\usepackage[francais]{babel}

\newcommand{\co}{\colon\thinspace}    
\newcommand{\fnote}[1]{\footnote{\small sharp1}}
\newcommand{\inv}{^{-1}}              

\newcommand{\N}{{\mathbb N}}
\newcommand{\Z}{{\mathbb Z}}
\newcommand{\R}{{\mathbb R}}
\newcommand{\Q}{{\mathbb Q}}

\newcommand{\T}{{\mathbb T}}

\newcommand{\spt}{\mbox{supp}}
\newcommand{\Azero}{\mathcal{A}_0}

\newtheorem{theorem}{Th\'eor\`eme}[chapter]
\newtheorem{proposition}[theorem]{Proposition}
\newtheorem{corollary}[theorem]{Corollaire}
\newtheorem{definition}[theorem]{D\'efinition}
\newtheorem{lemma}[theorem]{Lemme}
 
\newtheorem{conjecture}[theorem]{Conjecture} 
\newtheorem{question}[theorem]{Question}
\title{Syst\`emes lagrangiens et fonction $\beta$ de Mather}
\author{Daniel Massart}

\date{26 janvier  2011}

\begin{document}

\titlepage
\maketitle

\thispagestyle{empty}
\noindent{\Large \textbf{Remerciements et autres effusions}}
\vspace{1cm}

\noindent Albert Fathi m'a pris comme disciple en 1992, je me demande encore pourquoi. En 1997 il a bâti  en pierre de taille une théorie sur laquelle j'ai pu sculpter quelques gargouilles. En 2007 il m'a invité à bord du projet ANR KAM faible. Après quelques années dans le métier, je commence à mesurer la chance que j'ai eue en ces trois occasions, et l'ampleur de ma dette envers Albert.

L'ombre portée de John Mather s'étend si loin que j'ai bien du mal à écrire un article dont le titre ne contienne pas son nom. Il m'a fait, comme Gabriel Paternain,  un grand honneur en acceptant d'écrire un rapport sur ce mémoire. 

Patrick Bernard a passé des heures que j'imagine pénibles sur certains de mes articles, à en extirper impitoyablement erreurs et imprécision, me sauvant ainsi à plusieurs reprises de l'embarras éternel. De cela, et de s'être porté caution pour cette habilitation, je lui sais un gré infini.

Ivan Babenko, en plus d'être la vie et l'âme du séminaire Darboux, est la raison pour laquelle je ne me suis jamais découragé des mathématiques.
Je le remercie, ainsi que Marie-Claude Arnaud et Philippe Thieullen, d'avoir bien voulu faire partie du jury. 

La liste des mathématiciens que j'admire, et qui ont été bons pour moi au delà de mes mérites,
ne serait pas complète si elle ne mentionnait Pierre Arnoux et Victor Bangert. 

Un des bénéfices annexes des mathématiques est qu'elles m'ont permis de me lier d'amitié avec, par ordre approximativement chronologique, Frédéric Faure, Laurent Massoulié, le regretté Vincent Dumas, Jean-Marc Schlenker,  Ian Schindler, Philippe Narbel, T.R. Ramadas, Athoumane Niang, Mubariz Garaev. Tous ont contribué à ce travail à un moment ou un autre, en me remontant le Q.I., ou le moral.

Quoiqu'ils ne soient plus là pour me lire, je tiens à citer
\begin{itemize}
  \item Mr Yves Tanguy (1940-2007), instituteur à Locmaria-Plouzané, Finistère (plus largement  connu comme sonneur de bombarde). Il avait prédit que je serais mathématicien ou poète. Puisse-t-il y avoir encore des maîtres comme lui.
  \item mon père, Georges Massart (1939-2010). Il aurait aimé voir avancer ma carrière académique. Ce texte lui est dédié. 
 \end{itemize}

Enfin je n'oublie pas Jasmine, Etienne, et Emilie, qui me supportent tous les jours. 

\tableofcontents

\chapter{Introduction aux probl\`emes \'etudi\'es}
\textbf{Avertissement} : ce texte vise \`a pr\'esenter les r\'esultats de \cite{ijm}, \cite{soussol}, \cite{nonor}, \cite{vmb2}, \cite{AvsM}, \cite{Alfonso},   \cite{tworemarks}, \cite{codim1}, et \cite{completeproof}. Il ne constitue pas une introduction au sujet, ce pour quoi on consultera avec davantage de profit les r\'ef\'erences \cite{Bangert_twist}, \cite{Mather-Forni},   \cite{Fathi_bouquin}, et \cite{Mather_snowbird}. 

La bibliographie  recense en premier lieu les articles de ou cosign\'es par l'auteur, class\'es par ordre chronologique, et ensuite les articles cit\'es en r\'ef\'e\-rence, class\'es par ordre alphab\'etique sur le premier auteur. 
\section{Lagrangiens de Tonelli}

Les travaux expos\'es ici traitent de dyna\-mique des syst\`emes lagrangiens. Un tel syst\`eme est un mod\`ele de situation r\'egie par un principe d'\'economie, ou principe de moindre action. L'exemple canonique, et origine de la th\'eorie, est fourni par la m\'ecanique c\'eleste, o\`u le principe de moindre action fut \'enonc\'e par    Maupertuis, et formalis\'e ensuite par Lagrange. Un autre exemple, dont l'importance fut remarqu\'ee  par  Boltzmann (entre autres), est celui d'un r\'ecipient contenant un gaz parfait. Enfin les exemples abondent dans les activit\'es humaines, o\`u l'on cherche \`a minimiser un certain coût, exprim\'e en temps, en argent, ou en calories.

Formellement, on consid\`ere une vari\'et\'e diff\'erentiable $M$, appel\'ee espace de configurations, et une fonction $L$ sur le fibr\'e tangent de $M$, appel\'ee lagrangien.  
Traditionnellement le fibr\'e tangent \`a $M$ est appel\'e espace des phases. 
Dans le cas de la m\'ecanique c\'eleste, ou du gaz de Boltzmann, le lagrangien est la diff\'erence entre l'\'energie cin\'etique et l'\'energie potentielle. 
Nous serons souvent amen\'es \`a \'etudier le cas d'un syst\`eme soumis \`a une perturbation ext\'erieure p\'eriodique, telle que le passage d'une com\`ete dans le syst\`eme solaire. Le lagrangien est alors une fonction sur $TM \times \T$, o\`u $\T = \R / \Z$. 

Lorsque le lagrangien ne d\'epend pas du temps, il est dit autonome, et on omet le facteur $\T$.
L'\'evolution du syst\`eme est alors gouvern\'ee par l'\'equation d'Euler-Lagran\-ge, qui s'ex\-prime en cordonn\'ees locales $(x,v)$ sur $TM$, $t$ \'etant le temps, par 
\begin{equation}\label{Euler-Lagrange}
\frac{\partial}{\partial t} \frac{\partial L }{\partial v} (x,v,t) = \frac{\partial L }{\partial x} (x,v,t).
\end{equation}
Une d\'ecouverte importante de Lagrange  est que les projections dans $M\times \T$ des trajectoires (qui, elles, vivent dans $TM \times \T$) sont caract\'eris\'ees par le fait de minimiser localement, \`a extr\'emit\'es fix\'ees, l'int\'egrale du lagrangien, aussi appel\'ee int\'egrale d'action. Une courbe dans $M\times \T$ qui minimise localement, \`a extr\'emit\'es fix\'ees, l'int\'egrale  d'action, est appel\'ee extr\'emale. Une extr\'emale $C^1$ est donc la projection dans $M\times \T$ d'une trajectoire de l'\'equation d'Euler-Lagrange. Nous nous pla\c{c}ons en g\'en\'eral du point de vue des syt\`emes dynamiques, o\`u l'on \'etudie le comportement asymptotique des orbites. Pour cela il est commode de pouvoir garantir l'existence de solutions d\'efinies pour tous temps. Notons qu'en m\'ecanique c\'eleste, les solutions n'ont aucune raison d'\^etre d\'efinies pour tous temps, puisque les plan\`etes peuvent entrer en  collision. Sacrifiant le r\'ealisme \`a la simplicit\'e, nous ferons d\'esormais les hypoth\`eses suivantes, introduites par  Mather (\cite{Mather91}) :
\begin{itemize}
  \item (H1) $M$ est compacte, connexe, sans bord
  \item (H2) $L$ est de classe $C^2$
  \item (H3) stricte convexit\'e : la d\'eriv\'ee seconde $\partial^2 L / \partial v^2(x,v,t)$ est d\'efinie positive pour tout \'el\'ement $(x,v)$ de $TM$ et tout $t\in \T$
  \item (H4) surlin\'earit\'e :  pour tous $x$ dans $M$ et $t $ dans $\T$, $L(x,v,t) / \|v\| $ tend vers l'infini lorsque $ \|v\| $ tend vers l'infini, pour une m\'etrique riemannienne (quelconque, $M$ \'etant compacte)  sur $M$
  \item (H5) compl\'etude : le flot $\Phi_t$ sur $TM \times \T$ d\'efini 
  par l'\'equation d'Euler-Lagrange est complet, i.e. ses solutions maximales sont d\'efinies sur $\R$.
  \end{itemize}
Un lagrangien satisfaisant ces hypoth\`eses est dit de Tonelli.

Un dynamicien cherchera d\`es lors \`a  construire des orbites satisfaisant certaines propri\'et\'es globales, la plus simple \'etant d'\^etre p\'eriodique avec une  classe d'homotopie prescrite. Ce probl\`eme particulier est relativement facile et n'est mentionn\'e qu'\`a titre d'illustration. En revanche, la question suivante, pos\'ee par Arnold dans \cite{Arnold}, a inspir\'e les travaux fondateurs de Mather : 

\begin{question} Est-il vrai que pour un lagrangien autonome g\'en\'erique il existe une orbite dense dans presque tout  niveau d'\'energie ?
\end{question}
Dans le cas d'un lagrangien p\'eriodique  on peut se poser une  question analogue : 
\begin{question}\label{conjecture d'Arnold}
 Est-il vrai que pour un lagrangien g\'en\'erique il existe une orbite qui tend vers l'infini ?
\end{question}
Les hypoth\`eses de convexit\'e et surlin\'earit\'e (H3 et H4) nous invitent \`a employer des m\'ethodes variationnelles. La question est alors de trouver le bon espace dans lequel on cherche des extrema. Pour trouver des orbites p\'eriodiques d'homotopie prescrite, on cherchera naturellement \`a minimiser l'int\'egrale d'action sur l'espace des courbes $C^1$ dans la classe d'homotopie prescrite. Afin d'appliquer les m\'ethodes variationnelles \`a des probl\`emes plus g\'en\'eraux, Mather a propos\'e dans \cite{Mather91} la construction suivante.
\section{Th\'eorie de Mather}

 Notons $\mathcal{M}_{inv}$
l'espace des mesures bor\'eliennes de probabilit\'e  $\Phi_t$-invarian\-tes, \`a support compact dans  $TM\times \T$.
Mather a montr\'e que la fonction suivante, appel\'ee action du lagrangien sur les mesures,
	\[
	\begin{array}{rcl}
\mathcal{M}_{inv} & \longrightarrow & \R \\
\mu & \longmapsto & \int_{TM\times \T}	L d\mu
\end{array}
\]
est bien d\'efinie et admet un minimum. Une mesure r\'ealisant le minimum sera dite 
 $L$-minimisante.

Lorsque $M=\T$, par le   Th\'eor\`eme du Graphe de Mather, (\cite{Mather91}) une mesure invariante admet un nombre de rotation au m\^eme titre qu'une mesure invariante d'un hom\'eomorphisme du cercle. 

 Pour des syst\`emes ayant plus de degr\'es de libert\'e (c'est \`a dire que la dimension de $M$ est $>1$),  Mather a g\'en\'eralis\'e cette observation dans  \cite{Mather91} comme suit. 
Remarquons tout d'abord que si  $\omega$ est une 1-forme ferm\'ee sur  $M$ et $\mu \in \mathcal{M}_{inv}$ alors l'int\'egrale  $\int_{TM\times \T}	\omega d\mu$ est bien d\'efinie, et ne d\'epend que de la classe de cohomologie de 
 $\omega$. Par dualit\'e cela munit $\mu$ d'une classe d'homologie  $\left[\mu\right]$. La classe  $\left[\mu\right]$ est l'unique  $h \in H_1 (M,\R)$ tel que   
	\[
\langle h,\left[\omega \right]\rangle = \int_{TM\times \T}	\omega d\mu 
\]
pour toute 1-forme ferm\'ee $\omega$ sur  $M$. Il est prouv\'e dans  \cite{Mather91} que pour aucun  $h \in H_1 (M,\R)$, l'ensemble 
	\[ \mathcal{M}_{h,inv}:= \left\{\mu \in \mathcal{M}_{inv} \co \left[\mu\right]=h\right\}
\]
n'est vide. A nouveau  l'action du lagrangien sur ce sous-ensemble de mesures admet un minimum. La fonction qui \`a $h$ associe le minimum en question est appel\'ee la fonction  $\beta$ du  syst\`eme :
\[
	\begin{array}{rcl}
\beta \co H_1 (M,\R) & \longrightarrow & \R \\
h & \longmapsto & 
\min \left\{\int_{TM\times \T}	Ld\mu \co \mu \in \mathcal{M}_{h,inv} \right\}.
\end{array}
\] 
Une mesure $\mu$ telle que  $\left[\mu\right]=h$ et $\int L d\mu = \beta (h)$ est dite $(L,h)$-minimisante. 

On peut faire une construction analogue en cohomologie : si  $\omega$ est une 1-forme ferm\'ee sur   $M$, 
alors $L-\omega$ est un lagrangien satisfaisant les hypoth\`eses (H1-5), et de plus le flot d'Euler-Lagrange de  $L-\omega$ est le m\^eme que celui de  $L$. 
Le minimum sur  $\mathcal{M}_{inv}$ de $\int (L-\omega)d\mu$ ne d\'epend en fait que de la classe de cohomologie de  $\omega$, et d\'efinit une fonction sur $H^1(M,\R)$, dont l'oppos\'e est appel\'e fonction  $\alpha$ du syst\`eme :
\[
	\begin{array}{rcl}

\alpha \co H^1 (M,\R) & \longrightarrow & \R \\
c & \longmapsto & 
-\min \left\{\int_{TM\times \T}	(L-\omega)d\mu \co \mu \in  \mathcal{M}_{inv},\  \left[\omega\right]=c\right\}.
\end{array}
\]
Une mesure  $(L-\omega)$-minimisante est aussi dite $(L,\omega)$-minimisante, 
ou $(L,c)$-minimisante si $c$ est la classe de cohomologie de  $\omega$.

Mather a d\'emontr\'e que  $\alpha$ and $\beta$ sont  convexes, surlin\'eaires, et  duales de Fenchel, \`a savoir
\begin{eqnarray*}
\forall h \in H_1 (M,\R),\ \beta (h) &=& \sup_{c \in H^1 (M,\R)}\left( \langle c,h \rangle -\alpha (c) \right)\\
\forall c \in H^1 (M,\R),\ \alpha (c) &=& \sup_{h \in H_1 (M,\R)}\left( \langle c,h \rangle -\beta (h) \right).
\end{eqnarray*}
En particulier $\min \alpha = -\beta(0)$. 
Les principaux signes distinctifs d'une fonction convexe sont sa diff\'erentiabilit\'e et sa stricte (ou non) convexit\'e. Dans le cas des fonctions $\alpha$ et $\beta$, ils 
sont associ\'es \`a des propri\'et\'es dynamiques int\'eressantes, comme nous le verrons au paragraphe suivant. 

Un cas particulier digne d'int\'er\^et est celui o\`u le lagrangien est l'\'energie cin\'etique associ\'ee \`a une m\'etrique riemannienne ou finslerienne sur $M$. La sym\'etrie et l'homog\'en\'eit\'e du lagrangien se refl\`etent alors dans la fonction $\beta$, qui est, dans ce cas, le demi-carr\'e d'une norme, appel\'ee norme stable par Federer (cf.  \cite{Gromov-Lafontaine-Pansu}).

\section{Conjectures de Mather et Ma\~{n}\'e}
La strat\'egie de Mather pour attaquer la question \ref{conjecture d'Arnold}, aboutie en petites dimensions, et encore conjecturale en grandes,  consiste, sans entrer dans les d\'etails, \`a construire une sorte d'\'echelle infinie dont les barreaux sont des supports de mesures minimisantes. Les classes d'homologie de ces mesures minimisantes tendent vers l'infini, et l'orbite diffusante visite  tour \`a tour des voisinages de chacun des barreaux. Pour mener \`a bien cette construction, il importe que \begin{itemize}
  \item les mesures minimisantes ne soient pas trop nombreuses
  \item leurs supports ne soient pas trop gros.
  
\end{itemize}
Par exemple, quand $M=\T$, un tore invariant suffit \`a bloquer la diffusion puisqu'il s\'epare l'espace des phases en deux composantes connexes. Dans \cite{Mane95, Mane96, Mane97}, Ma\~n\'e a propos\'e une s\'erie de conjectures allant dans ce sens, pour des lagrangiens g\'en\'eriques. Pr\'ecisons d'abord ce que Ma\~{n}\'e entend par g\'en\'erique. Notons que si $L$ est un lagrangien de Tonelli sur $M$, et $f$ est une fonction (traditionnellement appel\'ee potentiel) $C^2$ sur $M\times \T$, alors $L+f$ est encore un lagrangien de Tonelli sur $M$. On dit alors qu'une propri\'et\'e vaut pour un lagrangien $C^k$-g\'en\'erique si, quel que soit un lagrangien de Tonelli $L$ de classe $C^k$, il existe une partie r\'esiduelle $\mathcal{O}(L)$ de $C^k (M\times \T)$ telle que pour tout $f \in \mathcal{O}(L)$, la propri\'et\'e vaut pour $L+f$. 
\begin{conjecture}[\cite{Mane97}]\label{Mane_forte}
Pour un lagrangien g\'en\'erique, il existe une uni\-que mesure minimisante, et cette mesure est port\'ee par une orbite p\'erio\-dique, ou un point fixe du flot d'Euler-Lagrange.
\end{conjecture}
Il existe une version faible de cette conjecture, o\`u on s'autorise \`a perturber le lagrangien par une 1-forme ferm\'ee en plus du potentiel :
\begin{conjecture}[\cite{Mane96}]\label{mane_faible}
Pour un lagrangien g\'en\'erique $L$, il existe un ouvert dense $U$ de $H^1(M,\R)$ tel que pour tout $c \in U$, il existe une   unique mesure $(L,c)$-minimisante, et cette mesure est port\'ee par une orbite p\'erio\-dique, ou un point fixe du flot d'Euler-Lagrange.
\end{conjecture}
A d\'efaut d'obtenir des orbites p\'eriodiques, on peut esp\'erer des mesures minimisantes ergodiques :
\begin{question}[\cite{Mane95}]\label{question_ergodicite}
Pour un lagrangien g\'en\'erique $L$, toute mesure minimisante est-elle ergodique ? 
\end{question}
Enfin on peut supposer que la coexistence de deux mesures minimisantes cohomologues est un ph\' enom\`ene exceptionnel (voir \cite{Bernard-Contreras} pour la meilleure r\' eponse \`a ce jour) : 
\begin{question}[\cite{Mane95}]\label{question_unicite_homologie}
Est-il vrai que pour un lagrangien g\'en\'erique $L$, pour toute classe d'homologie $h$, il existe  une unique mesure $(L,h)$-minimisante? 
\end{question}
\section{Diff\'erentiabilit\'e de $\beta$}
Dans \cite{Mather90} et \cite{Bangert94} un lien est mis en \'evidence entre petitesse des supports de mesures minimisantes et diff\'erentiabilit\'e de la fonction $\beta$. 
Le prototype de tous les th\'eor\`emes \`a ce sujet est 
\begin{theorem}[\cite{Mather90}, \cite{Bangert94}]\label{Mather1}
Si $M= \T$ alors $\beta$ est diff\'erentiable en toute classe  d'homologie irrationnelle. 
Elle est  diff\'erentiable en une classe d'homologie rationnelle $h$ si et seulement $\T$ est enti\`erement rempli d'orbites p\'eriodiques d'homologie $h$.
\end{theorem}
Puisque $H_1 ( \T,\R)=\R$, le mot  rationnelle se passe d'explication. Pour d'autres espaces de configurations il nous faut un peu de vocabulaire. 
Le quotient de  $H_1 (M,\Z)$ par sa torsion se plonge en un r\'eseau  $\Gamma$ dans $H_1 (M,\R)$. Une classe 
$h \in H_1 (M,\R)$ est dite enti\`ere si elle se trouve dans  $\Gamma$, et  rationnelle si  $nh \in \Gamma$ pour un entier non nul $n$. 
Un sous espace vectoriel de  $H_1 (M,\R)$ est dit entier s'il est engendr\'e par des classes enti\`eres.

Une fonction convexe admet un c\^one tangent en tout point. On dit qu'elle admet un sommet en $x$ si son c\^one tangent en $x$ ne contient aucune droite. Le th\'eor\`eme \ref{Mather1} sugg\`ere alors la question suivante : 
\begin{question} Est-il vrai que les sommets de $\beta$ ne se pr\'esentent qu'en des classes rationnelles ?
\end{question}

On aimerait donner un sens quantitatif \`a  l'irrationnalit\'e d'une classe d'homologie.
Le quotient $H_1 (M,\R)/\Gamma$ est un tore $\T^{b}$, o\`u $b$ est le premier nombre de  Betti de $M$. 
Pour $h$ dans $H_1 (M,\R)$, l'image de $\Z h$ dans $\T^{b}$ est un  sous-groupe de $\T^{b}$, 
donc son adh\'erence  $\mathcal{T}(h)$ est une r\'eunion finie de tores de m\^eme dimension. Cette  dimension est appel\'ee  irrationalit\'e de $h$, et not\'ee $I_{\Z}(h)$. 
Elle vaut  z\'ero si $h$ est rationnelle. Nous dirons qu'une classe $h$ est compl\`etement irrationnelle si son irrationalit\'e est maximale, \'egale \`a $b$.
Notons que  les irrationalit\'es de $h$ et  $nh$ sont \'egales, pour  $n\in \Z$, 
$n \neq 0$, puisque le quotient de $\mathcal{T}(h)$ 
par $\mathcal{T}(nh)$ est un  groupe fini.

Pour les lagrangiens autonome la notion d'irrationalit\'e pertinente est l\'eg\`erement diff\'erente. Soit $I_{\R}(h)$ la dimension de l'adh\'erence dans  $\T^{b}$ de  $\R h$ au lieu de  $\Z h$. Sur la relation entre les deux notions d'irrationalit\'e, voir l'appendice A.1 de \cite{vmb2}.  
Notons que la fonction $I_{\R}(h)$  est zero-homog\`ene, i.e. $I_{\R}(th)=I_{\R}(h)$ pour tous  $h \in H_1(M,\R)$ et $t \neq 0$.

On dit que $\beta$ est  diff\'erentiable dans $k$ directions en  $h$ si le c\^one tangent \`a $\beta$ en $h$ 
contient un espace affine de  dimension $k$. 
La question naturelle, que nous appellerons dans la suite le probl\`eme de diff\'erentiabilit\'e,  est alors 
 \begin{question} \label{question-diff de beta}
 La fonction $\beta$ est-elle  toujours  diff\'erentiable dans $k$ directions en une  classe d'homologie $k$-irrationnelle ?
 \end{question} 
Mather conjecture que la r\'eponse est oui pour les lagrangiens  $C^{\infty}$. 
La r\'eponse est oui pour tous les lagrangiens  $C^{2}$ quand  $M= \T$ d'apr\`es le Th\'eor\`eme \ref{Mather1}. 
\subsection{Int\'egrabilit\'e}
A l'origine du  probl\`eme de diff\'erentiabilit\'e se trouve la question suivante :  la non-diff\'erentiabilit\'e de $\beta$ ne se produit-elle  que lorsque les supports de mesures minimisantes sont petits ? A l'oppos\'e, on peut  s'interroger sur le rapport entre grosseur des supports et diff\'erentiabilit\'e de $\beta$. Plus pr\'ecisement 
\begin{question}\label{Burago-Ivanov} Si la fonction $\beta$ d'un lagrangien est $C^1$, le lagrangien est-il int\'egrable ?
\end{question}
Cette question est pos\'ee, dans le cadre des normes stables des m\'etriques de Finsler sur les tores, par Burago et Ivanov (cf. \cite{Burago-Ivanov}).
\section{Stricte convexit\'e de $\beta$ et ergodicit\'e des mesures minimisantes}
Tout d'abord un peu de vocabulaire : on appelle sous espace d'appui au graphe de $\beta$, un sous espace affine de $H_1(M,\R)\times \R$ qui rencontre le graphe de $\beta$, sans rencontrer l'\'epigraphe ouvert $\{ (h,y) \in H_1(M,\R)\times \R \co y > \beta(h) \}$. On appelle face de $\beta$ l'intersection du graphe de $\beta$ avec un sous espace d'appui. Une face est convexe. On dit que $\beta$ est strictement convexe en $h$ si la seule face de $\beta$ contenant $h$ est $\left\{h \right\}$.

Il est facile de voir que si $\beta$ poss\`ede une face non triviale, alors il existe une mesure minimisante non ergodique (combinaison convexe de mesures minimisantes pour des points extr\'emaux de la face). Inversement, si $\beta$ est strictement convexe en une classe d'homologie $h$, alors il existe une mesure minimisante ergodique d'homologie $h$. En effet toutes les composantes ergodiques d'une mesure $h$-minimisante quelconque ont pour homologie $h$, sans quoi il existerait une face de $\beta$ joignant les homologies des composantes. 

Ceci acquis, la r\'eponse \`a la question \ref{question_ergodicite} est contenue dans l'article \cite{Bangert90}, qui reprend une construction d'Hedlund (\cite{Hedlund}) : il existe une m\'etrique riemannienne sur le tore $\T^3$, dont la boule unit\'e de la norme stable est un octa\`edre. En particulier  chaque niveau d'\'energie non nul contient les supports d'exactement six mesures minimisantes ergodiques. On voit ais\'ement que cette situation persiste sous une perturbation du lagrangien par un potentiel, ce qui r\'epond \`a 
\ref{question_ergodicite}.

On peut alors se demander si une telle construction se g\'en\'eralise. Babenko et Balacheff ont montr\'e (\cite{BB06}) que sur toute vari\'et\'e de dimension $\geq 3$, il existe une m\'etrique dont la boule unit\'e de la norme stable est un poly\`edre, et de plus on dispose d'un choix assez large de poly\`edres. 

La situation est un peu plus compliqu\'ee en dimension deux. En effet, Bangert  a observ\'e que la boule unit\'e de la norme stable d'une m\'etrique sur $\T^2$ est toujours strictement convexe. Toutefois, si on consid\`ere des lagrangiens plus g\'en\'eraux que les m\'etriques, on dispose d'un degr\'e suppl\'ementaire de libert\'e puisque la fonction $\beta$ n'est en g\'en\'eral pas homog\`ene.  Dans \cite{CL99} Carneiro et Lopez ont mis en \'evidence un ph\'enom\`ene sp\'ecifiquement lagrangien : dans leur exemple, la fonction $\beta$ poss\`ede une face de la forme $\left[-h,h\right]$. Encore une fois cette situation est stable par perturbation.  

Les surfaces autres que le tore pr\'esentent des traits sp\'ecifiques, que nous d\'etaillons au chapitre  suivant.

\section{Th\'eorie d'Aubry-Mather en codimension un}
 La th\'eorie de Mather des lagrangiens de Tonelli (que nous appellerons th\' eorie de dimension un dans cette section, par opposition \`a la th\' eorie de codimension un) est une g\' en\' eralisation de la th\' eorie d'Aubry-Mather des twists maps de l'anneau, o\`u les solutions sont fonctions d'une variable, vue comme le temps, \`a valeurs dans l'espace de phase.  La th\' eorie d'Aubry-Mather en codimension un est une g\' en\' eralisation diff\' erente (et l\' eg\`erement plus ancienne puisque fond\' ee sur \cite{Moser} et \cite{Bangert unique}) : les solutions sont des fonctions \`a valeurs r\' eelles, et c'est la variable qui vit dans une vari\' et\' e de dimension quelconque, en l'occurrence  $\R^n$. Comme en th\' eorie de dimension un, les solutions sont caract\' eris\' ees par le fait de minimiser localement l'int\' egrale d'une fonction co\^ut, appel\' ee lagrangien. 
 
 Moser a mis en \' evidence l'existence de solutions dont le graphe reste \`a distance finie d'un hyperplan de $\R^n \times \R$. Cet hyperplan joue alors le r\^ole du nombre de rotation dans le cas des twist maps, ou de la classe d'homologie d'une mesure minimisante en th\' eorie de dimension un. On peut alors d\' efinir la fonction $\beta$ du syst\`eme : c'est la fonction qui \`a un hyperplan non vertical de $\R^n \times \R$ associe l'action moyenne minimale d'une solution dont le graphe est \`a distance finie de l'hyperplan en question. Par dualit\' e convexe on d\' efinit \' egalement la fonction $\alpha$.
 
Senn (\cite{Senn91})  a observ\' e que, comme dans le cas des twist maps, ou d'un lagrangien p\' eriodique en temps sur le cercle, ou d'un lagrangien autonome sur le tore de dimension deux, la fonction $\beta$ d'un syst\`eme lagrangien de codimension un est strictement convexe. Autrement dit $\alpha$ est $C^1$. Cela peut s'interpr\' eter de la fa\c{c}on suivante : dans tous les cas sus-cit\' es, les solutions tiennent suffisamment de place dans l'espace de configuration pour que deux d'entre elles, de classes d'homologie diff\' erentes, soient oblig\' ees de se couper. 

La diff\' erentiabilit\' e de $\beta$ en codimension un est \' etudi\' ee dans deux articles qui m\' eriteraient d'\^etre mieux connus, \cite{Senn95} et \cite{Bessi09}. 
En particulier la r\' eponse \`a la question \ref{question-diff de beta} est positive. Puisque une telle information est la premi\`ere \' etape vers la r\' esolution de la conjecture \ref{mane_faible}, on est alors tent\' e de r\' esoudre ladite conjecture dans le cadre de la th\'eorie d'Aubry-Mather en codimension un. Encore faut-il arriver \`a formuler la conjecture dans ce cadre. Il n'y a pas de dynamique en th\' eorie d'Aubry-Mather de codimension un, donc la notion de mesure invariante n'a pas de sens. Dans \cite{codim1}, \`a la suite de \cite{Bessi09}, lui-m\^eme inspir\' e par \cite{Bernard-Buffoni}, nous utilisons la notion de courant induit par une solution. 
Les probl\`emes s'\' enoncent comme suit. 
\begin{question}\label{question codim un}
Est-il vrai que pour un lagrangien  $L$ satisfaisant les hypoth\`eses de \cite{Moser}, g\' en\' erique en topologie $C^{\infty}$, 
\begin{itemize}
	\item 
	pour toute classe d'homologie $\rho$, les solutions $(L,\rho)$-minimisantes  induisent le m\^eme courant 
	\item pour toute classe de cohomologie $c$, les solutions $(L,c)$-minimisantes  induisent le m\^eme courant 
	\item il existe un ouvert dense $U$ de $\R^n$,  tel que pour tout  $c \in U$, vu comme une classe de cohomologie,  $\alpha'(c) \in \Q^n$ ?
	\end{itemize}
\end{question}
Le premier (resp. deuxi\`eme, resp. troisi\`eme) point  est une version de codimension un de la question \ref{question_unicite_homologie} (resp. \ref{question_ergodicite}, resp. \ref{mane_faible}).

\newpage
\chapter{R\'esultats de non-stricte convexit\'e}
Dans \cite{gafa} nous avons abord\'e le cas des normes stables des surfaces orientables. Un trait remarquable de ces derni\`eres est que leur homologie de dimension 1 est munie d'une structure symplectique canonique : la surface \'etant orient\'ee, on peut d\'efinir l'intersection alg\'ebrique de deux courbes param\'etr\'ees en position g\'en\'erale. Cette intersection ne d\'epend que des classes d'homologie des courbes. L'application $\mbox{Int} (.,.)$ ainsi d\'efinie sur  $H_1(M,\R) \times H_1(M,\R) $ est bilin\'eaire, antisym\'etrique, et non-d\'eg\'en\'er\'ee. 

Une cons\'equence du Graph Theorem de Mather (\cite{Mather91}) est que si deux classes d'homologie $h_1,h_2$ sont dans une face de $\beta$, alors $\mbox{Int} (h_1,h_2)=0$. En d'autres termes, les faces de la boule unit\'e de la norme stable sont contenues dans des sous-espaces isotropes de $H_1(M,\R)$. Ces derniers \'etant de dimension au plus la moiti\'e de celle de $H_1(M,\R)$, soit le genre de $M$, la boule unit\'e ne peut \^etre un poly\`edre. Ce fait avait d\'ej\`a \'et\'e observ\'e par M.J. Carneiro dans \cite{Carneiro}.

Notre contribution \`a ce sujet a consist\'e \`a montrer que la boule unit\'e n'est jamais strictement convexe d\`es que le genre de $M$ est $>1$, et plus pr\'ecis\'ement, pour toute classe d'homologie 1-irrationnelle $h$,  il existe  une face de la boule unit\'e, de dimension $\mbox{genre}(M)-1$, contenant $h$ (\cite{gafa}, th\'eor\`eme 7). En particulier la combinatoire des faces de la boule unit\'e est n\'ecessairement tr\`es compliqu\'ee. 

 Dans \cite{nonor}, en collaboration avec Florent Balacheff, nous avons \'etudi\'e la norme stable des surfaces non-orientables. Un trait qui diff\'erencie ces derni\`eres des orientables, et qui les rapproche des vari\'et\'es de dimension plus grande, est qu'on peut y tracer des courbes ferm\'ees disjointes dont les classes d'homologie engendrent $H_1(M,\R)$.

 Il est donc naturel que si $M$ est une surface non-orientable, il existe une 
  \newpage
\noindent m\'etrique sur $M$ dont la boule unit\'e de la norme stable est un poly\`edre :
 \begin{theorem}{(\cite{nonor}, th\'eor\`eme 1.3)} \label{nonor1}
 Soit  $M$ une  surface non-orien\-table ferm\'ee. Alors dans toute classe conforme il existe une m\'etrique dont la boule unit\'e de la norme stable est un poly\`edre.
\end{theorem} 
Cependant, contrairement \`a ce qui se passe en dimension plus grande, le nombre de sommets possible d'un tel poly\`edre est born\'e par $4b_1(M) -2$, o\`u $b_1(M) = \dim H_1 (M,\R)$. Cela tient au fait qu'on ne peut tracer sur $M$ plus de $2b_1(M) -1$ courbes ferm\'ees disjointes dont les classes d'homologie sont deux \`a deux non proportionnelles. 
 
 Par ailleurs la norme stable d'une  surfaces non-orientable $M$ munie d'une m\'etrique $g$ n'est jamais strictement convexe d\`es que $b_1(M) \geq 2$. Dans l'\'enonc\'e ci-dessous, $E$ d\'esigne la partie enti\`ere.  
  \begin{theorem}{(\cite{nonor}, th\'eor\`eme 1.1)}\label{nonor2}
 Soit $M$ une surface ferm\'ee non-orientable munie d'une m\'etrique riemannienne ou finslerienne. Alors toute courbe ferm\'ee minimisante est une composante connexe du support d'une mesure minimisante ayant au moins  $E((b_1(M)+1)/2)-1$ composantes ergodiques port\'ees par des courbes ferm\'ees.
\end{theorem} 
Cela s'explique par les consid\'erations suivantes.  Soient
 \begin{itemize}
  \item $p \co M_o \longrightarrow M$ le rev\^etement d'orientation de $M$
  \item $I$ l'involution canonique de $M_o$ induite par $p$
  \item $g_o$ la m\'etrique $p^* g$ sur $M_o$ 
  \item $E_1$ le sous espace de $H_1 (M_o,\R)$ form\'e des points fixes de $I_*$.
\end{itemize}
Alors $H_1(M,\R)$ s'identifie avec $E_1$, qui est un sous-espace lagrangien de $H_1(M_o,\R)$. La boule unit\'e de la norme stable de $(M,g)$ s'identifie avec l'intersection de $E_1$ avec la boule unit\'e de la norme stable de $(M_o,g_o)$.

Dans \cite{AvsM}, nous avons g\'en\'eralis\'e le th\'eor\`eme 7  de \cite{gafa} au cas de lagrangiens quelconques. Il ne s'agit pas de g\'en\'eraliser pour le plaisir, ce r\'esultat est n\'ecessaire \`a la preuve du th\'eor\`eme 1.3 de \cite{AvsM}, lui-m\^eme utilis\'e dans la preuve de la conjecture de Ma\~n\'e (voir Annexe C), ainsi que dans \cite{Alfonso}. Dans le cas d'une m\'etrique, le lagrangien est homog\`ene, toute l'information sur la fonction $\beta$ est contenue dans la boule unit\'e.  En particulier les faces de $\beta$ sont contenues dans des ensembles de niveaux de $\beta$. Dans le cas g\'en\'eral, le th\'eor\`eme est un peu plus compliqu\'e \`a \'enoncer.
 
\begin{definition}
Soit  $M$ une vari\'et\'e ferm\'ee et soit  $L$ un lagrangien de  Tonelli sur  $TM\times \T$.
Pour  $h \in H_1(M,\R)\setminus \left\{0\right\}$, nous d\'efinissons la face radiale maximale de $\beta$ contenant $h$, not\'ee $R_h$,  comme \'etant le plus grand sous-ensemble de la demi-droite  $\left\{th \co t \in \left[0,+\infty \right[ \right\}$ contenant  $h$ (pas n\'eces\-sairement dans son int\'erieur relatif) en restriction auquel  $\beta$ est affine.
\end{definition}

 La possibilit\'e de faces radiales est la diff\'erence la plus flagrante  entre les fonctions $\beta$ des m\'etriques riemanniennes ou finsleriennes (\cite{gafa}, \cite{nonor}) et celles des lagrangians g\'en\'eraux. Un exemple de face radiale est exhib\'e dans  \cite{CL99}. Nous d\'efinissons l'ensemble de  Mather  $\tilde{\mathcal{M}}(R_h)$  comme la r\'eunion de tous les  supports de mesures  $th$-minimisantes, pour $th \in R_h$.
 
 On note $V(R_h)$ le sous-espace vectoriel engendr\'e, dans  $H_1(M,\R)$ par les classes d'homologies de toutes les composantes ergodiques de mesures $th$-minimisantes, pour $th \in R_h$ (ce sont les mesures dont le support est contenu dans $\tilde{\mathcal{M}}(R_h)$). Si $V$ est un sous-espace vectoriel de $H_1(M,\R)$, $M$ \'etant une surface orient\'ee, on note $V^{\perp}$ l'orthogonal de $V$ pour la forme symplectique sur $H_1(M,\R)$ induite par l'intersection alg\'ebrique des courbes.  

Enfin, on dit qu'une classe de cohomologie est non-singuli\`ere si son ensemble d'Aubry (voir le sous-chapitre \ref{Aubry_definition} pour la d\'efinition de l'ensemble d'Aubry) ne contient aucun point fixe du flot d'Euler-Lagrange. Si $h$ est une classe d'homologie, la transform\'ee de Legendre $\mathcal{L}(h)$ est une face de la fonction $\alpha$, elle a donc un ensemble d'Aubry (voir le sous-chapitre \ref{faces_de_alpha} pour le rapport entre ensembles d'Aubry et faces de $\alpha$). On dit qu'une classe d'homologie est non-singuli\`ere si l'ensemble d'Aubry de sa transform\'ee de Legendre ne contient aucun point fixe du flot d'Euler-Lagrange. 
La raison de ces pr\'ecautions pour \'eviter les points fixes est que nos d\'emonstrations reposent sur le fait que, en dimension deux, les courbes ferm\'ees sont localement s\'eparantes.
\begin{theorem}[(\cite{AvsM}, th\'eor\`eme C.4)]\label{AvsM_faces}
Soient
\begin{itemize}
	\item $M$ une surface ferm\'ee, orient\'ee
	\item $L$ un lagrangien de  Tonelli sur $M$
	\item $h_0$ un \'el\'ement 1-irrationnel, non-singulier de  $H_1(M,\R)$
	\item $V_0 := V(R_{h_0})$
	\item $h \in V^{\perp}_{0}$.
\end{itemize}
Alors il existe $t(h_0,h) \in \R$ et $s(h_0,h) >0$ tels que le  segment \newline
$\left[h_0,t(h_0,h)h_0 + s(h_0,h)h\right]$ est contenu dans une  face de $\beta$.
\end{theorem}

La d\'emonstration de ce th\'eor\`eme repose sur les lemmes et propositions suivants. 
Rappelons que si $h$ est une classe d'homologie, 
$$\mathcal{L}(h) = \left\{c \in H^1 (M,\R) \co \alpha(c)+ \beta(h) = \left\langle c,h \right\rangle\right\}$$
 est le sous-diff\'erentiel \`a $\beta$ en $h$, et $\mathcal{L}(h)$ est une face de $\alpha$.
\begin{lemma}[\cite{AvsM}, lemme 2.3] \label{key_o}
Soient \begin{itemize}
  \item $M$ une surface ferm\'ee, orient\'ee
  \item $L$ un lagrangien de  Tonelli autonome sur $M$
  \item $\gamma_0$ une extr\'emale ferm\'ee minimisante de  $L$, telle que  $\gamma_0$ n'est pas un  point   fixe
  \item $h_0$ la classe d'homologie de la mesure minimisante  port\'ee par  
  $(\gamma_0,\dot\gamma_0)$
  \item $c$ une  classe de cohomologie dans  $\mathcal{L}(h_0)$.
    \end{itemize}
 Alors il existe un voisinage $V_0$ de $(\gamma_0,\dot\gamma_0)$ dans $TM$ telle que pour toute extr\'emale   $\gamma$ telle que  $(\gamma,\dot\gamma) \subset \tilde{\mathcal{A}}(c)$, si  $(\gamma,\dot\gamma)$ entre dans  (resp. sort de ) $V_0$  alors  $\gamma$ est  
 \begin{itemize}
  \item ou bien une  extr\'emale ferm\'ee homotope \`a  $\gamma_0$
  \item ou bien  positivement (resp. negativement ) asymptote \`a une  extr\'emale ferm\'ee homotope \`a $\gamma_0$.
 \end{itemize}
  \end{lemma}
  Rappelons que l'int\'erieur relatif d'un convexe d'un espace de Banach est son int\'erieur dans le sous-espace affine qu'il engendre. Il est d\'emontr\'e dans \cite{gafa}, voir aussi \cite{AvsM}, que si deux faces $F_1,F_2$ contiennent un point donn\'e $x$ dans leurs int\'erieurs relatifs, alors il existe une face contenant $F_1$ et $F_2$, et contenant $x$ dans son int\'erieur relatif. Il y a par cons\'equent un sens \`a parler de la plus grande face contenant un point donn\'e dans son int\'erieur.
\begin{lemma}[\cite{AvsM}, lemme C.4]\label{AvsM_C4}
Soit  $F$ une face de  $\beta$.  Supposons que  $F$ contienne une classe d'homologie  $1$-irrationnelle  $h_0$ dans son  int\'erieur relatif. Alors pour tout  $h \in F$, pour toute mesure $h$-minimisante  $\mu$, le support de $\mu$ consiste en orbites  p\'eriodiques, et/ou points fixes.
\end{lemma}
Le lemma \ref{AvsM_C4} est cons\'equence de la 
\begin{proposition}[\cite{CMP}, Proposition 2.1  ]\label{rational}
Soit $M$ une  surface fer\-m\'ee, pas n\'ecessairement orientable, et soit  $L$ un lagrangien de  Tonelli sur $M$.
Si  $h$ est une classe d'homologie 1-irrationnelle et $\mu$ est une mesure $h$-minimiante, alors le  support de $\mu$ consiste en orbites  p\'eriodiques, ou (non exclusif) points fixes.
\end{proposition}
Cette proposition a une histoire un peu compliqu\'ee : elle m'a \'et\'e expliqu\'ee par A. Fathi, qui en attribuait l'argument \`a Haefliger. Rassur\'e par cette auguste caution, j'en ai \'ecrit dans \cite{gafa} une vague esquisse de preuve. Les auteurs de \cite{CMP}, ennemis de l'\`a-peu-pr\`es, ont r\'edig\'e une preuve compl\`ete. Ignorant ce fait, j'ai \'ecrit avec F. Balacheff une autre preuve compl\`ete  dans \cite{nonor}.  
\section{Un probl\`eme de comptage}
Soit $(M,g)$ une surface ferm\'ee munie d'une m\'etrique riemannienne ou finslerienne. Nous appellerons 
\begin{itemize}
	\item $\mathcal{S}(g)$ l'ensemble des longueurs de g\'eod\'esiques ferm\'ees qui sont dans le support d'une mesure $h$-minimisante, $h$ \'etant une classe d'homologie 1-irrationnelle, compt\'ees avec multiplicit\'es, c'est \`a dire que si deux g\'eod\'esiques non homologues ont m\^eme longueur, celle-ci est compt\'ee deux fois 
 \item $N(T)$ le nombre d'\'el\'ements de $\mathcal{S}(g)$ qui sont $\leq T$.
\end{itemize}
Le probl\`eme qui nous int\'eresse ici est de trouver des estimations asymptotiques de $N(T)$. Voyons quelques exemples pour commencer.
\subsection{Tores}
Supposons que $M=\T^2$. Si $h$ est une classe d'homologie enti\`ere  non nulle et $\mu$ une mesure $h$-minimisante, alors $\mu$ est port\'ee par des g\'eod\'esiques ferm\'ees d'apr\`es la proposition \ref{rational}. De plus ces g\'eod\'esiques ferm\'ees sont deux \`a deux disjointes en vertu du Graph Theorem de Mather, et non s\'eparantes puisque $h \neq 0$. Cela implique qu'elles sont deux \`a deux homotopes. Puisqu'elles minimisent la longueur dans leur classe d'homologie, elles ont toutes m\^eme longueur. Ainsi, \`a toute classe d'homologie enti\`ere on peut associer un \'el\'ement de $\mathcal{S}(g)$. R\'eciproquement tout \'el\'ement de de $\mathcal{S}(g)$ est la norme stable d'une classe d'homologie enti\`ere. Donc $N(T)$ est \'egal au nombre $N_1 (T)$ de classes d'homologie enti\`eres dont la norme stable est $\leq T$. Par un th\'eor\`eme classique de Minkowski, 
	\[\lim_{T \rightarrow \infty}\frac{N_1 (T)}{T^2} = \mathcal{V}(g) 
\]
o\`u $\mathcal{V}(g)$ d\'esigne le volume de la boule unit\'e de la norme stable, calcul\'e au moyen de la forme volume sur $H_1 (M,\R)$ pour laquelle le r\'eseau $\Gamma$ est de d\'eterminant $1$.  On est alors tent\'e de conjecturer que pour une surface de premier nombre de Betti $b$, on a 
\[\lim_{T \rightarrow \infty}\frac{N_1 (T)}{T^b} = \mathcal{V}(g).
\]

\subsection{Surfaces non-orientables}
Le pr\'esent paragraphe  a pour objet de refroidir notre enthousiasme. Commen\c{c}ons par observer que si $M$ n'est pas le tore, et si $\gamma$ est une g\'eod\'esique ferm\'ee minimisant la longueur dans sa classe d'homologie $h$, alors la boule unit\'e de la norme stable poss\`ede un sommet en $h / \|h\|$ (voir \cite{gafa}, th\'eor\`eme 8). 
D'apr\`es le th\'eor\`eme \ref{nonor1}, si $M$ n'est pas orientable, il existe une m\'etrique $g$ sur $M$ dont la boule unit\'e de la norme stable est un poly\`edre. En particulier elle a un nombre fini de sommets, donc $\mathcal{S}(g)$ est fini. Nous fermerons d\'esormais les yeux sur ce navrant \'etat de fait, et supposerons $M$ orientable. Rappelons (voir \cite{gafa}) que dans ce cas, si $M$ est de genre strictement plus grand que $1$, la boule unit\'e de la norme stable a toujours une infinit\'e de sommets. 
\subsection{Girafes}
Supposons que $M$ est de genre $k$, et poss\`ede $k$ g\'eod\'esiques ferm\'ees s\'eparantes $\delta_1, \ldots, \delta_k$,  ayant la propri\'et\'e suivante : chaque $\delta_i$ poss\`ede un voisinage tubulaire $V_i$ \`a bord lisse, tel que tout arc $C^1$ qui traverse $V_i$ d'un bord \`a l'autre est au moins aussi long que n'importe lequel des deux bords. Remarquons  que les m\'etriques ayant cette propri\'et\'e, ou girafes,  forment une partie d'int\'erieur non vide, et d'adh\'erence non compacte, de l'espace des m\'etriques sur $M$. Alors aucune g\'eod\'esique minimisante ne traverse un des $V_i$. En effet, si une g\'eod\'esique $\gamma$ traverse $V_i$, en rempla\c{c}ant les arcs de $\gamma$ qui sont contenus dans $V_i$ par des segments du bord de $V_i$, on obtient un cycle rectifiable, homologue \`a $\gamma$, et plus court. Les $\delta_i$ d\'ecoupent $M$ en $k$ tores priv\'es d'un disque $T_1, \ldots, T_k$. On a alors
\begin{itemize}
	\item $H_1(M,\R) = \oplus^{k}_{i=1} H_1 (T_i, \R)$
	\item $H_1 (T_i, \R)$ est un sous-espace symplectique de dimension deux de \\
	$H_1(M,\R)$, pour tout $i$, et pour la forme symplectique induite en homologie par l'intersection alg\'ebrique des courbes
	\item la boule unit\'e $\mathcal{B}$ de la norme stable de $M$ est l'enveloppe convexe des boules unit\'es $\mathcal{B}_i$ des normes stables des $T_i$, munis de la m\'etrique induite par celle de $M$ (en particulier, les seuls sommets de $\mathcal{B}$ sont ceux des $\mathcal{B}_i$)

   \item $N(T)\approx \sum^{k}_{i=1} \mathcal{V}_i T^2$, o\`u $\mathcal{V}_i$ d\'esigne le volume de $\mathcal{B}_i$, calcul\'e au moyen de la forme symplectique d'intersection.
 \end{itemize} 
 
 \subsection{Fraises}
 Dans ce paragraphe on aborde le cas g\'en\'eral. 
 Nous quittons ici la terre ferme pour nous lancer dans la sp\'eculation. Supposons pour simplifier que $M$ est de genre deux. Soit $\gamma_0$ une g\'eod\'esique ferm\'ee minimisant la longueur dans sa classe d'homologie, et soit $h_0$ la classe d'homologie $\left[\gamma_0 \right]/l_g(\gamma_0)$, o\`u $l_g$ d\'esigne la longueur pour la m\'etrique $g$. Ainsi $h_0$ est un sommet de la boule unit\'e de la norme stable. Alors (\cite{gafa}, th\'eor\`eme 7, voir aussi \cite{nonor}, th\'eor\`eme 6.1) pour toute classe d'homologie $h_1$ telle que $\mbox{Int} (h_0,h_1)=0$, $\mbox{Int}$ d\'esignant la forme symplectique d'intersection,  il existe $s=s(h_1,h_0)$ tel que le segment de droite joignant $h_0$ \`a $\left(h_0 + sh_1\right)/\|h_0 +s h_1\|$ est contenu dans le bord de la boule unit\'e de la norme stable. Nous appelons fraise (au sens vestimentaire et non mara\^\i cher) de $h_0$ la partie suivante de la boule unit\'e :
	\[ \mathcal{F}(h_0) := \left\{ h_1 \in H_1(M,\R) \co \forall t \in \left[0,1\right],\  \|th_0 + (1-t)h_1 \|=1 \right\}.
\]
Nous dirons que deux fraises sont \'equivalentes si leurs bords coïncident. 
Par exemple dans le cas o\`u $(M,g)$ est une girafe de genre deux, il y a deux classes d'\'equivalences de fraises, correspondant aux  intersections de la sph\`ere unit\'e avec les sous-espaces symplectiques $H_1 (T_i, \R)$. Dans le cas d'une surface de genre $k$ quelconque, on associe une fraise \`a toute classe d'homologie enti\`ere poss\'edant un cycle minimisant \`a $k-1$ composantes connexes. Ainsi  une girafe de genre $k$ poss\`ede exactement $k$ classes d'\'equivalences de fraises, correspondant aux  intersections de la sph\`ere unit\'e avec les sous-espaces symplectiques $H_1 (T_i, \R)$.

En g\'en\'eral on peut montrer que  $\mathcal{F}(h_0)$ est hom\'eomorphe, de fa\c{c}on bi-lipschitzienne, \`a un disque. Il y a donc un sens \`a int\'egrer la forme symplectique d'intersection sur une fraise.  On peut montrer que deux fraises \'equivalentes ont la m\^eme aire symplectique. Appelons $\mathcal{F}(g)$ la somme, dans $\left]0, +\infty\right]$, des aires symplectiques de toutes les classes d'\'equivalence de fraises. 

Ne risquant rien d'autre que le ridicule, nous proposons la
\begin{conjecture}
Pour toute surface riemannienne ou finslerienne $(M,g)$, on a
	\[\lim_{T \rightarrow \infty}\frac{N_1 (T)}{T^2} = \mathcal{F}(g).
\]
\end{conjecture}
Naturellement cet \'enonc\'e n'a d'int\'er\^et que si $\mathcal{F}(g) < \infty$, ce que nous nous gardons bien de conjecturer. 
\chapter{Diff\'erentiabilit\'e de $\beta$ et ensembles d'Aubry}
\section{Sous-solutions de l'\'equation d'Hamilton-Jacobi}\label{Aubry_definition}
Commen\c{c}ons par quelques rudiments de th\'eorie KAM faible de Fathi (voir \cite{Fathi_bouquin} pour plus de d\'etails). Adoptons, pour le moment, le point de vue hamiltonien. Le hamiltonien $H$ associ\'e \`a $L$ est
$$
\begin{array}{rcl}
H \co T^* M \times \T & \longrightarrow & \R \\
(x,p,t) & \longmapsto & \sup_{v \in T_x M} \left\langle p,v \right\rangle - L(x,v,t).
\end{array}
$$
La convexit\'e et la surlin\'earit\'e de $L$ garantissent que $H$ est bien d\'efini, strictement convexe dans les fibres de $T^* M$, et surlin\'eaire.  La transform\'ee de Legendre
$$
\begin{array}{rcl}
\mathcal{L} \co T M \times \T & \longrightarrow & T^* M \times \T \\
(x,v,t) & \longmapsto & (x, \frac{\partial L }{\partial v}(x,v,t), t)
\end{array}
$$
est un diff\'eomorphisme $C^1$ qui conjugue le flot d'Euler-Lagrange de $L$ avec le flot hamiltonien de $H$. Une m\'ethode classique pour trouver des ensembles invariants non triviaux du flot consiste \`a r\'esoudre l'\'equation de Hamilton-Jacobi
\begin{equation}\label{Hamilton-Jacobi}
\frac{\partial u}{\partial t}(x,t) + H(x,d_x u, t) = \mbox{constante}
\end{equation}
o\`u l'inconnue $u$ est une fonction $C^1$ de $M \times \T$ vers $\R$. Faute de solutions, en g\'en\'eral, on cherche des sous-solutions, i.e. des fonctions  $u \in C^1 (M \times \T, \R) $ telles que
\begin{equation}\label{Hamilton-Jacobi_ineq}
\frac{\partial u}{\partial t}(x,t) + H(x,d_x u, t) \leq \mbox{constante}.
\end{equation}
Soit 
$$
I := \left\{ c \in \R \co \exists u \in C^1 (M \times \T, \R), \;  \forall x \in M,\; \frac{\partial u}{\partial t}(x,t) + H(x,d_x u, t) \leq c \right\}.
$$
L'ensemble $I$ est non-vide puisque n'importe quelle fonction est sous-solution pour un $c$ assez grand. D'autre part toute sous-solution pour $c$ est sous-solution pour $c' \geq c$, donc $I$ est un intervalle infini \`a droite. Par compacit\'e de $M$, convexit\'e et surlin\'earit\'e de $H$,  $H$ est born\'e inf\'erieurement, ce qui, dans le cas o\`u $H$ est autonome, entra\^\i ne que $I$ est born\'e inf\'erieurement. Dans le cas non-autonome, c'est un lemme dû \`a Albert Fathi (voir \cite{soussol}). Appelons valeur critique, et notons $\alpha_H$, la borne inf\'erieure. Observons tout d'abord que la notation $\alpha_H$ n'est pas fortuite : en effet, Fathi a  montr\'e que $\alpha_H = \alpha(0)$ o\`u $\alpha(0)$ d\'esigne la valeur en la classe de cohomologie nulle de la fonction $\alpha$ du syst\`eme. 

La question est alors de savoir si $I$ est ferm\'e. La r\'eponse est apport\'ee par le th\'eor\`eme suivant, dû \`a Fathi et Siconolfi dans le cas autonome et \`a l'auteur dans le cas p\'eriodique.
\begin{theorem}[\cite{FS}, \cite{soussol}]\label{FS04}
Il existe $u \in C^1 (M \times \T, \R)$ telle que 
$$ \forall x \in M,\  \frac{\partial u}{\partial t}(x,t) + H(x,d_x u, t) \leq \alpha_H.
$$
\end {theorem}
On appelle sous-solution critique une telle fonction $u$. Notons qu'il existe n\'ecessairement une partie $\mathcal{E}_u$ de $M \times \T$ o\`u l'\'egalit\'e est r\'ealis\'ee dans l'in\'egalit\'e de Hamilton-Jacobi, sans quoi on contredirait la minimalit\'e de $\alpha_H$.  Ceci nous am\`ene \`a la d\'efinition suivante :
\begin{definition}[Fathi]
On appelle ensemble d'Aubry de $H$, et on note $\mathcal{A}(H)$, l'intersection sur toutes les sous-solutions critiques $u$ des $\mathcal{E}_u$.
\end{definition}
Voici pourquoi il valait la peine de s'int\'eresser aux sous-solutions, et de   d\'efinir l'ensemble d'Aubry :
\begin{proposition}[Fathi]
(1) L'ensemble $\mathcal{A}(H)$ est non vide et ferm\'e.

\noindent (2) Quelles que soient $u_1$ et $u_2$ deux sous-solutions critiques, et quel que soit $(x,t) \in \mathcal{A}(H)$, on a 
$$
\frac{\partial u_1}{\partial t}(x,t) = \frac{\partial u_2}{\partial t}(x,t), \mbox{ et } d_x u_1 = d_x u_2.
$$ 

\noindent (3) La partie 
$$
\left\{ (x, d_x u, t) \co (x,t) \in \mathcal{A}(H), \ u \mbox{ sous-solution critique } \right\}
$$
 de $T^* M \times \T$ est invariante par le flot hamiltonien.
\end{proposition}
Dans la suite on appellera ensemble d'Aubry du lagrangien $L$ associ\'e \`a $H$ par dualit\'e de Fenchel, et on notera $\tilde{\mathcal{A}}(L)$, l'image de 
$$
\left\{ (x, d_x u, t) \co (x,t) \in \mathcal{A}(H), \ u \mbox{ sous-solution critique } \right\}
$$
 dans $TM \times \T$ par la transform\'ee de Legendre.

En fait le th\'eor\`eme \ref{FS04} dit un peu plus : il affirme l'existence de sous-solutions strictes, c'est \`a dire que l'in\'egalit\'e de Hamilton-Jacobi est stricte hors de l'ensemble d'Aubry. Une cons\'equence int\'eressante pour nous est la 
\begin{proposition}\label{H+W}
Il existe une fonction $W \in C^2 (M\times \T, \R)$, telle que $W \geq 0$, $W^{-1}(0) = \mathcal{A}(H)$, et l'ensemble d'Aubry $\mathcal{A}(H+W)$
est \'egal \`a $\mathcal{A}(H)$.
\end{proposition}
Cette proposition a un sens puisque $H+W$ est le hamiltonien associ\'e au lagrangien de Tonelli $L-W$. Voyons tout d'abord comment d\'eduire la proposition \ref{H+W} du th\'eor\`eme \ref{FS04}. Soit $u$ une sous-solution stricte $C^1$, et soit $W$ une fonction $C^2$ telle que
	\[ \forall x,t,\ \frac{\partial u}{\partial t}(x,t) + H(x,d_x u, t)+W(x,t) \leq \alpha_H.
\]
Alors $u$ est sous-solution pour le hamiltonien de Tonelli $H+W$, donc  $\alpha_{H+W}\leq \alpha_H$.
D'autre part toute sous-solution pour $H+W$ est \'evidemment sous-solution pour $H$, donc $\alpha_{H+W}\geq \alpha_H$. Par cons\'equent $\alpha_{H+W}=  \alpha_H$, et $\mathcal{A}(H+W)\supset \mathcal{A}(H)$. Mais puisque $u$ est stricte en dehors de $\mathcal{A}(H)$, on a $\mathcal{A}(H+W)\subset \mathcal{A}(H)$, ce qui d\'emontre la proposition.

Dans \cite{soussol} on fait le chemin inverse : on d\'emontre d'abord la proposition \ref{H+W}, \`a l'aide de la d\'efinition lagrangienne de l'ensemble d'Aubry, et ensuite on en d\'eduit le th\'eor\`eme \ref{FS04}.

Rappelons que si $\omega$ est une 1-forme ferm\'ee sur $M$, alors $L-\omega$ est un lagrangien de Tonelli, dont le flot d'Euler-Lagrange est le m\^eme que celui de celui de $L$. De plus il est \'el\'ementaire de v\'erifier que l'ensemble d'Aubry de $L-\omega$ ne d\'epend que de la classe de cohomologie $c$ de $\omega$. Si $H_{\omega}$ d\'esigne le hamiltonien associ\'e \`a $L-\omega$, on a $\alpha_{H_{\omega}}= \alpha (c)$.
On est alors amen\'e \`a se demander comment varie l'ensemble d'Aubry en fonction de $c$. 
\section{Faces de $\alpha$}\label{faces_de_alpha}
Voici le rapport entre  la discussion pr\'ec\'edente et la dif\-f\'erentiabilit\'e de $\beta$. L'id\'ee de base nous a \'et\'e inspir\'ee par la lecture de \cite{Mather93}. Observons que par dualit\'e convexe, il est \'equivalent d'\'etudier la diff\'erentiabilit\'e de $\beta$ et la stricte convexit\'e de $\alpha$. Les propositions suivantes relient les faces de $\alpha$ avec les variations de l'ensemble d'Aubry vu comme fonction de la classe de cohomologie. La notion d'int\'erieur relatif d'une face joue un rôle central. 
\begin{proposition}[\cite{ijm}, \cite{vmb2}]\label{Aubry_constant_faces}
Soient
\begin{itemize}
	\item $M$ une vari\'et\'e ferm\'ee
	\item $L$ un lagrangien de Tonelli sur $TM\times \T$
	\item $F$ une face de $\alpha$. 
\end{itemize}
Si $(c_0, \alpha(c_0))$ est un point de l'int\'erieur relatif de $F$, et si $(c, \alpha(c))$ est un point quelconque de $F$, alors l'ensemble d'Aubry de $c$ contient celui de $c_0$. En particulier l'ensemble d'Aubry, vu comme fonction de la classe de cohomologie, est constant dans l'int\'erieur relatif de $F$.  
\end{proposition}
Nous appellerons d\'esormais ensemble d'Aubry d'une face $F$ de $\alpha$ l'ensemble d'Aubry commun \`a tous les points de l'int\'erieur relatif de $F$.
La proposition \ref{Aubry_constant_faces} admet une r\'eciproque partielle : 
\begin{proposition}[\cite{ijm}]\label{Aubry_constant_faces_reciproque}
Soient
\begin{itemize}
	\item $M$ une vari\'et\'e ferm\'ee
	\item $L$ un lagrangien de Tonelli sur $TM\times \T$
	\item $c$ et $c'$ deux classes de cohomologie dont les ensembles d'Aubry se rencontrent.
\end{itemize}
Alors il existe une face de $\alpha$ contenant $(c,\alpha(c))$ et $(c',\alpha(c'))$.
\end{proposition}
 Des  propositions \ref{H+W} et \ref{Aubry_constant_faces_reciproque} nous pouvons d\'eduire le 
 \begin{lemma}
Soient
\begin{itemize}
	\item $M$ une vari\'et\'e ferm\'ee
	\item $L$ un lagrangien de Tonelli sur $TM\times \T$
	\item $\omega$ une 1-forme ferm\'ee sur $M \times \T$ dont le support ne rencontre pas l'ensemble d'Aubry de $L$
	\item $(c,\tau)$ la classe de cohomologie de $\omega$ dans $H^1(M,\R) \times \R$.
\end{itemize} 
Alors il existe $\epsilon > 0$, tel que l'ensemble d'Aubry de $L\pm \epsilon \omega$ est \'egal \`a celui de $L$, et $\alpha$ est affine en restriction au segment $\left[-\epsilon c, \epsilon c\right]$.
\end{lemma}
Ce lemme peut se reformuler de fa\c{c}on plus conceptuelle.   

\begin{definition} 
Soit $E_0$ l'ensemble des $(c,\tau) \in H^1 (M\times \T,\R)=H^1 (M,\R)\times H^1 ( \T,\R)$ tels qu'il existe une 1-forme ferm\'ee lisse  $\omega$ sur $M\times \T$ avec $[\omega] = (c,\tau)$ et $\spt (\omega)\cap \tilde{\mathcal{A}}_0 = \emptyset$. 
\end{definition}
\begin{definition}
Soit  $G_0$ l'ensemble des  $(c,\tau) \in H^1 (M\times \T,\R)=H^1 (M,\R)\times H^1 ( \T,\R)$ tels qu'il existe une 1-forme ferm\'ee continue  $\omega$ sur $M\times \T$ avec $[\omega] = (c,\tau)$ et 
	\[\omega ((x,t),(v,\tau))=0 \;\forall (x,t) \in \Azero \subset M\times \T, \; \forall (v,\tau)\in 
	T_{(x,t)}M\times \T.
\]
\end{definition}
\begin{definition}
Soit $V_0$ le sous-espace vectoriel sous-jacent \`a l'espace affine engendr\'e par la plus grande face (qui existe d'apr\`es \cite{gafa}) de $\alpha$ contenant zero dans son int\'erieur. 
\end{definition}
On a alors le th\'eor\`eme suivant, d\'emontr\'e dans \cite{ijm} pour le cas autonome, et dans \cite{soussol} pour le cas p\'eriodique.
\begin{theorem}
Soient
\begin{itemize}
	\item $M$ une vari\'et\'e ferm\'ee
	\item $L$ un lagrangien de Tonelli sur $TM\times \T$.
	\end{itemize}
On a les inclusions suivantes : $E_0 \subset V_0 \subset G_0$.
\end{theorem}
L'\'egalit\'e $E_0= V_0$ entra\^\i ne une r\'eponse affirmative au probl\`eme de dif\-f\'erentiabilit\'e. En effet $E_0$ est un sous-espace entier (c'est \`a dire un sous-espace engendr\'e par des classes enti\`eres) de $H^1(M,\R)\times \R$, puisque $E_0$ s'identifie \`a la cohomologie d'un ouvert de $M$, laquelle est engendr\'ee par des classes enti\`eres par le th\'eor\`eme des coefficients universels. D'autre part, on a la
\begin{proposition}[\cite{vmb2}, proposition 20]\label{vmb2_prop20}
Soit  $L$ un lagrangien de \\
 Tonelli sur une vari\'et\'e ferm\'ee $M$.
Supposons que pour toute classe de cohomologie $c \in H^1 (M,\R)$,  $c$, $\tilde{V}_c$ est un sous-espace entier de 
$H^1 (M\times \T^1,\R)$. Soit  $h$ une classe d'homologie  $k$-irrationnelle. 
Alors $\beta$ est diff\'erentiable en  $h$ en au moins $k$ directions.
\end{proposition}
D\'esormais le but du jeu sera donc pour nous de trouver des conditions suffisantes \`a l'\'egalit\'e $E_0= V_0$. Notons que ce ne saurait \^etre le cas 
g\'en\'eral, puisqu'un contre-exemple est construit dans \cite{BIK}. Comme souvent dans ce sujet, les probl\`emes se posent lorsque l'ensemble d'Aubry est de grande dimension de Hausdorff, et  de forme compliqu\'ee (par opposition \`a rectifiable).
\section{Application jacobienne}
Notre outil pour chercher des cas d'\'egalit\'e $E_0=V_0$, est l'application, dite jacobienne, introduite dans \cite{vmb2},  et d\'efinie comme suit.
Soit $\omega$ une 1-forme ferm\'ee sur $M\times \T$ telle que l'ensemble d'Aubry de $L-\omega$ est \'egal \`a celui $\mathcal{A}(L)$  de $L$. Appelons $p \co \overline{M}\times \R \longrightarrow M\times \T$ le rev\^etement ab\'elien de $M\times \T$, c'est \`a dire le plus petit rev\^etement dans lequel le rappel de toute 1-forme ferm\'ee sur $M\times \T$  est exact. Soient  
\begin{itemize}
	\item $F \co \overline{M}\times \R \longrightarrow \R$ une primitive du rappel de $\omega$ \`a $\overline{M}\times \R$
	\item $x_0$ un point quelconque de $\mathcal{A}(L)$
	\item $u_0$ et $u_1$ des sous-solutions de l'\'equation de Hamilton-Jacobi associ\'ees \`a $L$ et $L-\omega$, respectivement.
\end{itemize}
On d\'efinit alors l'application jacobienne associ\'ee \`a $\omega$ par
	\[ \begin{array}{rcl}
	\phi_{\omega} \co \overline{M}\times \R & \longrightarrow & \R  \\
	x & \longmapsto & \left( F+u_1\circ p -   u_0\circ p  \right)(x,t) - 
     \left( F+u_1\circ p -   u_0\circ p  \right)(x_0,0).
     \end{array}.
\]
Une propri\'et\'e importante, pour nous, de cette application, est donn\'ee par la 
\begin{proposition}[\cite{vmb2}, Proposition 6]\label{Holder}
L'application $\phi_{\omega}$ satisfait une condition de    H\H{o}lder d'ordre deux sur $p\inv (\Azero) $.
 \end{proposition}
 
\section{Sommets de $\beta$}
Un cas particulier du probl\`eme de diff\'erentiabilit\'e, est la question suivante, d\'ej\`a cit\'ee dans l'introduction : \textit{ est-il vrai que les sommets de $\beta$ ne se pr\'esentent qu'en des classes rationnelles ?}
Le th\'eor\`eme 1 de \cite{ijm} entra\^\i ne une r\'eponse affirmative quand la dimension de $M$ est deux, mais on peut dire un peu plus par une m\'ethode sp\'ecifique : 
\begin{theorem}[\cite{vmb2}, th\'eor\`eme 3]\label{vmb2_vertices}
Si la fonction $\beta$ d'un lagrangian $L$ sur une vari\'et\'e  $M$ a un sommet en une classe d'homologie $h$, et $p$ est l'irrationalit\'e de $h$, alors  $2p < \dim M$. 
\end{theorem}

L'id\'ee de la d\'emonstration est la suivante. Supposons, sans perte de g\'en\'eralit\'e, que la transform\'ee de Legendre $\mathcal{L}(h)$ de $h$ est la plus grande face de $\alpha$ contenant z\'ero dans son int\'erieur. Si on a un sommet de $\beta$ en $h$, alors on peut trouver une base $c_1, \ldots c_b$ de $H^1 (M,\R)$, form\'ee de classes de cohomologie enti\`eres, et des $\epsilon_1, \ldots \epsilon_b >0$  tels que $\epsilon_1 c_1, \ldots \epsilon_b  c_b \in \mathcal{L}(h)$. Pour chaque $i=1, \ldots b$ on choisit une 1-forme ferm\'ee $\omega_i$ repr\'esentant $\epsilon_i c_i$, et on forme l'application 
	\[ \begin{array}{rcl}
\Phi \co	\overline{M}& \longrightarrow & \R^b \\
	x & \longmapsto & \left(\phi_{\omega_i} \right)_{i=1, \ldots b}
	\end{array}
\]
o\`u $\phi_{\omega_i}$ est l'application jacobienne associ\'ee \`a $\omega_i$ comme au paragraphe pr\'ec\'edent. Puisque les 1-formes ferm\'ees $\epsilon_i \inv \omega_i$, $i=1, \ldots b$, sont enti\`eres, cette application passe au quotient en une application $\Phi$ de $M$ dans un tore $\T^b$.  Si $\gamma \co \R \longrightarrow TM \times \T$ est une orbite du flot d'Euler-Lagrange de $L$,  contenue dans l'ensemble d'Aubry $\mathcal{A}(L)$, alors $\gamma(\Z) \subset M \times \left\{0\right\} \subset M \times \T$, et en identifiant $M \times \left\{0\right\}$ \`a $M$, on peut parler de $\Phi\left(\gamma(\Z)\right)$. Alors $\Phi\left(\gamma(\Z)\right)$ est dense dans un sous-tore de $\T^b$, dont la dimension est \'egale \`a l'irrationnalit\'e de $h$. Puisque $\Phi$ est 2-h\H{o}lderienne le long de  $\mathcal{A}(L)$, on conclut par le lemme de Ferry que $2p < \dim M$.

Comme les pr\'ec\'edents, ce th\'eor\`eme admet un analogue dans le cas autonome. 
Observons que pour un lagrangien autonome, d'apr\`es \cite{Carneiro}, $\beta$ est toujours diff\'erentiable dans la direction radiale, c'est \`a dire que, pour $h \in H_1 (M,\R) \setminus \left\{0\right\}$ donn\'e, la fonction  
\[
\begin{array}{rcl}
\left]0,+\infty \right[ & \longrightarrow & \left]0,+\infty \right[ \\ 
 t & \mapsto & \beta (th)
 \end{array}
 \]
  est $C^{1}$. Donc la fonction  $\beta$ d'un lagrangien autonome  ne peut avoir de sommets qu'en z\'ero. 
Le minimum de diff\'erentiabilit\'e possible dans le cas autonome est donc atteint quand $\beta$ n'est d\'erivable en une classe $h \neq 0$ que dans la direction radiale, auquel cas le c\^one tangent \`a $\beta$ en $h$ ne contient aucun plan. Notre r\'esultat, qui am\'eliore le th\'eor\`eme 1 de \cite{BIK}, est le 
\begin{theorem}\label{vmb2_vertices_autonome}
Si la fonction  $\beta$ d'un lagrangien autonome $L$ sur une vari\'et\'e ferm\'ee $M$ n'est diff\'erentiable en  $h \neq 0$ que dans la direction radiale,    
alors  $2I_{\R}(h)-1 < \dim M$. 
En particulier, si la  dimension de $M$ est trois, alors $I_{\R}(h)\leq 1$. 
\end{theorem}
Ce th\'eor\`eme est d\'emontr\'e dans \cite{vmb2} sous l'hypoth\`ese qu'il n'y a pas de points fixes dans l'ensemble d'Aubry. Cette hypoth\`ese est  superflue ;  une preuve est donn\'ee en annexe \ref{vmb2+}.
\section{Questions}
\subsection{Flots d'Anosov}
Que peut-on dire si le flot d'Euler Lagrange est d'Anosov ?
\subsection{Petites dimensions, g\'en\'ericit\'e }
On constate une analogie formelle entre le probl\`eme de diff\'erentiabilit\'e et le probl\`eme de la dimension de Hausdorff de l'ensemble d'Aubry quotient. Ce dernier probl\`eme \'etant r\'esolu en petites dimensions dans  \cite{FFR} ( \`a la suite de \cite{Mather02} qui r\'esout le probl\`eme de sa totale discontinuit\'e), et dans \cite{Bernard-Contreras}, pour un lagrangien g\'en\'erique en toutes dimensions,  on aimerait un r\'esultat analogue pour 
le  probl\`eme de diff\'erentiabilit\'e.

\subsection{Adh\'erence de courbes sur les tores}\label{question_courbe}
En voulant g\'en\'eraliser les th\'eor\`emes \ref{vmb2_vertices} et \ref{vmb2_vertices_autonome}, on peut chercher \`a affaiblir l'hypoth\`ese de non-diff\'erentiabilit\'e : par exemple dans le cas p\'eriodique en temps, au lieu de supposer que $\beta$ a un sommet en $h$, on suppose seulement que le cône tangent \`a $\beta$ en $h$ ne contient aucun plan. L'image par $\Phi$ d'une orbite contenue dans l'ensemble d'Aubry est alors contenue, non plus dans une droite de $\T^b$, mais dans un plan. On est donc amen\'e \`a consid\'erer la question suivante : soit $\gamma$ une courbe de $\R^n$, trac\'ee sur un sous-espace totalement irrationnel de $\R^n$ (c'est \`a dire un sous-espace ne contenant aucun point \`a coordonn\'ees enti\`eres, si ce n'est l'origine), et admettant une direction asymptotique (le quotient $\gamma(t) / t$ tend vers une limite non nulle quand $t$ tend vers l'infini). Que peut-on dire sur la dimension de Hausdorff de l'image de $\gamma$ dans $\R^n / \Z^n$ ? Sauf dans le cas facile o\`u $\gamma$ est \`a distance born\'ee d'une droite de $\R^n$ (voir annexe \ref{reponse_courbe}), cette question demeure myst\'erieuse.

\chapter{Ensembles d'Aubry et de Mather pour un lagrangien autonome  en dimension deux}

Si on dispose d'une r\'eponse affirmative au probl\`eme de diff\'erentiabilit\'e, et si on veut conna\^\i tre le c\^one tangent \`a $\beta$ en une classe d'homologie $h$, il reste \`a obtenir suffisament d'information sur l'ensemble d'Aubry pour pouvoir d\'eterminer $E_0$. Rappelons que si $h$ est une classe d'homologie, $\mathcal{L}(h) = \left\{c \in H^1 (M,\R) \co \alpha(c)+ \beta(h) = \left\langle c,h \right\rangle\right\}$ est le sous-diff\'erentiel \`a $\beta$ en $h$, et $\mathcal{L}(h)$ est une face de $\alpha$. On note $\mathcal{A}\mathcal{L}(h)$ l'ensemble d'Aubry commun \`a tous les points de l'int\'erieur relatif de $\mathcal{L}(h)$. De m\^eme on note  $\mathcal{M}\mathcal{L}(h)$ l'ensemble de Mather commun \`a tous les points de l'int\'erieur relatif de $\mathcal{L}(h)$, c'est \`a dire la r\'eunion de tous les supports de mesures $c$-minimisantes, pour $c$ dans l'int\'erieur relatif de  $\mathcal{L}(h)$. En g\'en\'eral $\mathcal{M}\mathcal{L}(h)\subset \mathcal{A}\mathcal{L}(h)$, et $\mathcal{M}\mathcal{L}(h)\setminus \mathcal{A}\mathcal{L}(h)$ est fait d'orbites homoclines \`a $\mathcal{M}\mathcal{L}(h)$. Notre th\'eor\`eme suivant permet d'exclure les homoclines dans certains cas particuliers. 

On dit qu'une classe d'homologie $h$ est singuli\`ere si $\mathcal{M}\mathcal{L}(h)$ contient un point fixe du flot d'Euler-Lagrange. Naturellement cela n'a de sens que pour un lagrangien autonome. Notre r\'esultat s'appuie fortement sur le fait que les courbes ferm\'ees sont localement s\'eparantes en dimension deux, c'est pourquoi il nous faut en exclure les classes singuli\`eres.
\begin{theorem}\label{AvsM}
Soient
\begin{itemize}
	\item $M$ une  surface ferm\'ee
	\item $L$ un lagrangien de  Tonelli autonome sur  $M$
	\item $h$ une classe d'homologie 1-irrationnelle, non singuli\`ere.
\end{itemize}
Alors $\mathcal{AL}(h)=\mathcal{ML}(h)$, et $\mathcal{AL}(h)$ est une r\'eunion d'orbites p\'eriodiques.
\end{theorem}
Ce r\'esultat poss\`ede un analogue dans le cas des lagrangiens p\'eriodiques sur le cercle, \'enonc\'e dans \cite{Mather93}, p. 1376, et d\'emontr\'e dans \cite{Mather_snowbird}, section 3. Soit $L$ un lagrangien de Tonelli sur $T\T \times \T$, et soit $h$ une classe d'homologie rationnelle de $\T$. Alors $\mathcal{L}(h)$ est un intervalle $\left[c^-,c^+ \right]$, et pour tout $c$ dans l'int\'erieur relatif de $\left[c^-,c^+ \right]$ (\`a savoir $\left]c^-,c^+ \right[$ si $c^+ \neq c^-$, et $\left\{c^+\right\}$ si $c^+ = c^-$), on a $\mathcal{A}(c)=\mathcal{M}(c)$, et $\mathcal{A}(c)$ est une r\'eunion d'orbites p\'eriodiques.

Outre le fait, \'enonc\'e par Poincar\'e,  que les orbites p\'eriodiques sont la br\`eche par laquelle nous pouvons esp\'erer investir la forteresse des syst\`emes dynamiques, nous avons des raisons pr\'ecises pour chercher des orbites p\'erio\-diques minimisantes, \`a savoir : 
\begin{itemize}
	\item les conjectures de Ma\~n\'e
	\item les r\'esultats de  \cite{A03} et \cite{AIPS05}
	\item le r\'esultat de \cite{Bernard07}, qui affirme l'existence de sous-solutions de l'\'equa\-tion de Hamilton-Jacobi aussi r\'eguli\`eres que le Lagrangien, pourvu que l'ensemble d'Aubry consiste en une r\'eunion finie d'orbites p\'eriodiques hyperboliques.
\end{itemize}
Pour les deux derniers points il est essentiel que l'ensemble d'Aubry ne contienne que des orbites p\'eriodiques, auquel cas il est \'egal \`a l'ensemble de Mather. 

Une premi\`ere observation est qu'en dimension deux il est relativement facile de trouver des orbites p\'eriodiques minimisantes : la proposition \ref{rational} nous assure que si $h$ est une classe d'homologie 1-irrationnelle, toute mesure $h$-minimisante est port\'ee par des orbites p\'erio\-diques ou points fixes. Cela ne suffit pas \`a garantir que $\mathcal{ML}(h)$ soit form\'e d'orbites p\'eriodiques. En effet, si $c \in \mathcal{L}(h)$, on pourrait a priori imaginer qu'il  existe dans $\mathcal{L}(c)$ une classe d'homologie $h'$ $k$-irrationnelle avec $k>1$, et une mesure $h'$-minimisante qui ne soit pas port\'ee par des orbites p\'eriodiques. En ce cas $h$ et $h'$ sont contenues dans $\mathcal{L}(c)$ qui est une face de $\beta$. Cela nous oblige \`a analyser avec pr\'ecision les faces de $\beta$, et justifie le temps pass\'e \`a d\'emontrer le th\'eor\`eme \ref{AvsM_faces}. 

\begin{definition}
Si $L$ est un lagrangien de Tonelli (autonome ou p\'eri\-odique) sur une vari\'et\'e ferm\'ee $M$, et si $F$ est une face de $\alpha$, on pose 
	\[\mathcal{F}(F) := \bigcap_{c \in F}\mathcal{L}(c).
\]
\end{definition}
Si $h$ est une classe d'homologie de $M$, $\mathcal{L}(h)$ est une face de $\alpha$, on note $\mathcal{FL}(h) := \mathcal{F}\left(\mathcal{L}(h)\right)$ pour \'eviter la surcharge. Alors $\mathcal{FL}(h)$ est une face de $\beta$ (\cite{AvsM}, lemme A.1). Il est clair que $h$ appartient \`a  $\mathcal{FL}(h)$, mais pas forc\'ement \`a son int\'erieur relatif. Si on pouvait montrer que $h$ est contenu dans l'int\'erieur relatif de $\mathcal{FL}(h)$,  par le lemme \ref{AvsM_C4} on obtiendrait que $\mathcal{ML}(h)$ est form\'e d'orbites p\'eriodiques. Dans le cas d'un lagrangien p\'eriodique sur le cercle, la fonction $\beta$ est strictement convexe, donc  $\mathcal{FL}(h)=\left\{h \right\}$. Pour un lagrangien autonome, la situation est compliqu\'ee par la possibilit\'e de faces radiales. 

Soit  $h$ une classe d'homologie 1-irrationnelle. Rappelons que $R_h$ est la plus grande face radiale de $\beta$ contenant $h$ (pas n\'ecessairement dans son int\'erieur). 
Alors pour tout  $t$ tel que  $th \in R_h$, $th$ est aussi  1-irrationnelle. 
S'il est faux en g\'en\'eral que $h$ soit contenue dans l'int\'erieur relatif de $\mathcal{FL}(h)$, le lemme suivant est suffisant pour nos projets : 
\begin{lemma}[\cite{AvsM}, lemme 3.1]\label{nonsingular_1}
Soit  $L$ un lagrangien autonome  sur une  surface ferm\'ee $M$.
Soit  $h$ une classe d'homologie  $1$-irrationnelle, non singuli\`ere. Alors l'int\'e\-rieur relatif de  $R_h$ est contenu dans l'int\'erieur relatif de   $\mathcal{FL}(h)$.
\end{lemma}
Nous savons \`a pr\'esent que $\mathcal{ML}(h)$ est form\'e d'orbites p\'eriodiques, pour peu que $h$ soit 1-irrationnelle et non singuli\`ere. Il reste \`a prouver que $\mathcal{AL}(h)=\mathcal{ML}(h)$. L'id\'ee est la suivante. Soit $\gamma$ une orbite p\'eriodique minimisante. Supposons pour simplifier que $M$ est orientable, $\gamma$ est non-s\'eparante, et l'ensemble d'Aubry ne contient pas d'orbites p\'eriodiques qui s'accumulent sur $\gamma$. Soient 
\begin{itemize}
	\item $h := \left[\gamma \right]$
	\item $h_0$ une classe d'homologie dont l'intersection alg\'ebrique avec $h$ vaut un
	\item $h^{\pm}_{n} := nh \pm h_0$, pour tout entier $n >0$
	\item des mesures  $h^{\pm}_{n}$-minimisantes  $\mu^{\pm}_{n}$
	\item des classes de cohomologie $c^{\pm}_{n}$ telles que  $\mu^{\pm}_{n}$  soit $c^{\pm}_{n}$-minimisante pour tout $n$
	\item des valeurs d'adh\'erence $c^{\pm}$ des suites $c^{\pm}_{n}$.
\end{itemize}
Alors en passant \`a la limite dans l'in\'egalit\'e de Fenchel on voit que $c^{\pm} \in \mathcal{L}(h)$. D'autre part, en vertu de la semi-continuit\'e de l'ensemble d'Aubry, valable en dimension deux d'apr\`es \cite{Bernard_Conley}, toute valeur d'adh\'erence, au sens de la topologie de Hausdorff, des suites des supports des mesures $\mu^{\pm}_{n}$
est contenue dans $\mathcal{A}(c^{\pm})$. On en d\'eduit que les ensembles d'Aubry $\mathcal{A}(c^{\pm})$ contiennent des orbites homoclines \`a $\gamma$ de tous les côt\'es, et dans tous les sens. Puisque $\mathcal{A}\mathcal{L}(h)$ est contenu \`a la fois dans $\mathcal{A}(c^{+})$ et $\mathcal{A}(c^{-})$, et v\'erifie le Graph Theorem de Mather, cela exclut que $\mathcal{A}\mathcal{L}(h)$ contienne des orbites homoclines \`a $\gamma$. L'ensemble d'Aubry \'etant constitu\'e de l'ensemble de Mather et d'orbites homoclines \`a ce dernier, on obtient ainsi que $\gamma$ est isol\'ee dans  $\mathcal{A}\mathcal{L}(h)$. Plus g\'en\'eralement, on a la 
\begin{proposition}[\cite{AvsM}, Corollaire 2.5] \label{AvsM_corollaire_local}
Soit  $c$ dans $H^1 (M,\R)$ telle que l'ensemble de  Mather  $\tilde{\mathcal{M}}(c)$ consiste en orbites p\'eriodiques qui ne sont pas des points fixes,  $\gamma_i, i \in I$. Soit $h$ un barycentre \`a coefficients strictement positifs des classes d'homologie des $\gamma_i, i \in I$. Alors 
	\[\mathcal{AL}(h)= \mathcal{ML}(h)= \mathcal{M}(c).
\]
\end{proposition}
Ceci ach\`eve la d\'emonstration du th\'eor\`eme \ref{AvsM}. 
\section{Question}
Les exemples expos\'es dans \cite{AvsM} montrent que l'hypoth\`ese de non-singularit\'e est n\'ecessaire au th\'eor\`eme \ref{AvsM}. Plus pr\'ecis\'ement, elle est n\'ecessaire \`a la partie de l'\'enonc\'e qui concerne l'ensemble d'Aubry : on peut construire des lagrangiens dont l'ensemble d'Aubry contient des orbites homoclines, dont on ne peut se d\'ebarasser en d\'epla\c{c}ant la classe de cohomologie \`a l'int\'erieur d'une face. Toutefois on peut se demander si l'hypoth\`ese de non-singularit\'e est n\'ecessaire \`a la partie de l'\'enonc\'e qui concerne l'ensemble de Mather : serait-il possible de construire un lagrangien sur un tore de dimension deux, dont l'ensemble de Mather contiendrait un point fixe et une lamination sans feuille ferm\'ee, tel que la fonction $\beta$ serait diff\'erentiable en z\'ero ? 
\section{Cons\'equence du th\'eor\`eme \ref{AvsM} sur la diff\'erentiabilit\'e de $\beta$}
Dans le cadre de la th\'eorie d'Aubry-Mather de codimension un, Senn a obtenu  la superbe formule suivante :
\begin{theorem}[\cite{Senn95}]
Soient 
\begin{itemize}
  \item $L$ un lagrangien satisfaisant les hypoth\`eses de Moser (\cite{Moser})
  \item $h \in \R^{n}$, on voit alors $(h,1)$  comme une classe d'homologie de dimension $n$ sur le tore $\T^{n+1}$
  \item $\mathcal{V}(h,1) := \mbox{Vect}\{ c_1- c_2   \co   c_1, c_2 \in \mathcal{L}(h,1) \}$
  \item $\mbox{rat}(h,1) := \mbox{Vect} (h,1)^{\perp} \bigcap \Z^{n+1}$, autrement dit, $\mbox{rat}(h,1)$ est la partie enti\`ere de l'orthogonal euclidien de $(h,1)$, on le voit comme une partie de la cohomologie de dimension $n$ du tore $\T^{n+1}$.
 \end{itemize} 
Alors 
$$\mathcal{V}(h,1)= \mbox{rat}(h,1), $$
\`a moins que les solutions minimisantes d'homologie $(h,1)$ ne feuillettent $\T^{n+1}$, auquel cas $\mathcal{V}(h,1)=\{0\}$ (i.e. $\beta$ est diff\'erentiable en $(h,1)$).
\end{theorem}
En d'autres termes on a la dichotomie suivante : ou bien les solutions minimisantes feuillettent, auquel cas $\beta$ est diff\'erentiable ind\'ependamment de la rationalit\'e de $h$, ou bien elles ne feuillettent pas, auquel cas la diff\'erentiabilit\'e de $\beta$ est enti\`erement d\'etermin\'ee par la rationalit\'e de $h$. 

Dans le cadre de dimension un, on aimerait bien exprimer la diff\'erentiabilit\'e de $\beta$ sous une forme aussi compacte. Il n'y a aucun espoir d'\'etendre telle qu'elle la formule de Senn : si $M$ est une surface de genre $>1$, munie d'un lagrangien de Tonelli $L$, et si $h$ est une classe d'homologie rationnelle de $M$, la diff\'erentiabilit\'e de $\beta$ en $h$ ne d\'epend pas seulement de la rationalit\'e de $h$, mais aussi du lagrangien $L$, via le nombre de composantes connexes des mesures $(L,h)$-minimisantes. Dans \cite{Alfonso} nous obtenons la formule suivante : 
\begin{theorem}\label{betadiff, thm}
Soient
\begin{itemize}
	\item $M$ une  surface ferm\'ee
	\item $L$ un lagrangien de   Tonelli autonome sur  $M$
	\item $h_0$ une classe d'homologie  1-irrationelle, non singuli\`ere  de $M$
	\item $(\gamma_i, \dot\gamma_i)_{i \in I}$ les orbites p\'eriodiques contenues dans les supports des mesures 
	 $(L,th)$-minimisantes, pour tout  $th$ in $R_{h_0}$
	\item $c_0$ une classe de cohomologie contenue dans l'int\'erieur relatif de $\mathcal{L}(h_0)$
	\item $\mathcal{V}(h_0)$ le sous-espace de  $H^1(M,\R)$ engendr\'e par les diff\'erences $c_1 - c_2$, o\`u $c_1,c_2 \in \partial \beta (h_0)$
\end{itemize}
Alors on a :
\begin{itemize}
  \item ou bien 
  $$
  \mathcal{V}(h_0) = \partial \alpha (c_0)^{\perp}= \bigcap_{i \in J } h_i ^{\perp}
  $$
  o\`u l'orthogonalit\'e s'entend au sens de la dualit\'e entre  $H_1(M,\R)$ et $H^1(M,\R)$
  \item ou bien $M= \T^2$   et les courbes ferm\'ees $\gamma_i $, $i \in I$, feuillettent  $M$;  dans ce cas $ \mathcal{V}(h_0) = \{0\}$. 
\end{itemize}
\end{theorem}

\section{Perspective : orbites homoclines}
Supposons donn\'e un lagrangien autonome sur une surface ferm\'ee.
Nous venons de voir que si une classe de cohomologie se trouve dans l'int\'erieur relatif d'une face rationnelle de $\alpha$, alors son ensemble d'Aubry consiste en orbites p\'eriodiques. R\'eciproquement, si l'ensemble d'Aubry d'une classe de cohomologie $c$ consiste en orbites p\'eriodiques, alors il est facile de voir que $c$ se trouve dans l'int\'erieur relatif d'une face de $\alpha$. Maintenant, consid\'erons une face rationnelle quelconque de $\alpha$. Sur son bord se trouvent des points extr\'emaux. L'ensemble d'Aubry d'un tel point extr\'emal ne peut \^etre enti\`erement constitu\'e d'orbites p\'eriodiques, d'apr\`es les consid\'erations ci-dessus. Par cons\'equent, il contient des orbites homoclines aux orbites p\'eriodiques qui forment l'ensemble d'Aubry de l'int\'erieur de la face. Nous avons ainsi une m\'ethode nouvelle pour trouver des orbites homoclines, qui de surcro\^\i t sont contenues dans des ensembles d'Aubry. 
\chapter{Int\'egrabilit\'e}
Dans \cite{Alfonso} nous abordons, en collaboration avec A. Sorrentino,  la question \ref{Burago-Ivanov}. Notre r\'eponse, partielle, met en jeu la notion de $C^0$-int\'egrabilit\'e. 
\begin{definition}[\cite{Arnaud}]\label{C0integrability}
Un hamiltonian de Tonelli $H: T^*M \longrightarrow \R$ est dit $C^0$-integrable, s'il existe un feuilletage de  $T^* M$ dont les feuilles sont des graphes lagrangiens,  invariant, lipschitziens, et dans chaque classe de  cohomologie se trouve une feuille.
\end{definition}
Notre r\'esultat principal est le 
\begin{theorem}[\cite{Alfonso}, th\'eor\`eme 1]\label{Alfonso_maintheo1}
Soit   $L \co T \T^2 \longrightarrow \R$ un lagrangien de  Tonelli autonome sur le tore de dimension deux. 
Alors la fonction $\beta$ est $C^1$ si et seulement si le syst\`eme est $C^0$-int\'egrable.
\end{theorem}
Observons que le tore est la seule surface qui admette un lagrangien $C^0$-int\'egrable (\cite{Alfonso}, Proposition 4). 
D'autre part   si $M$ est la sph\`ere de dimension deux, le plan projectif r\'eel, ou la bouteille de Klein, la fonction $\beta$ est $C^1$ quel que soit le lagrangien (trivialement pour les deux premi\`eres, d'apr\`es \cite{Carneiro} pour la derni\`ere). Si $M$ n'est ni  la sph\`ere de dimension deux, le plan projectif r\'eel,  la bouteille de Klein, ou le tore de dimension deux, la fonction $\beta$ n'est jamais $C^1$ quel que soit le lagrangien.
C'est une cons\'equence de la proposition suivante, qui g\'en\'eralise un r\'esultat de Bangert pour les flots g\'eod\'esiques sur le tore de dimension deux, et  de Mather pour les twist maps. 
\begin{proposition}[\cite{Alfonso}, proposition 1]\label{diff1irrat} Soit $L$ un lagrangien de  Tonelli sur une  surface ferm\'ee $M$ et soit   $h \in H_1(M, \R)$ une classe d'homologie  $1$-irrationnelle non singuli\`ere. Alors , $\beta$ est diff\'erentiable en $h$  si et seulement si $\mathcal{M} (R_h)=M$. En ce  cas, $\mathcal{M} (R_h)=M$ est feuillet\'e par des extr\'emales ferm\'ees. 
\end{proposition}
Pour d\'emontrer le th\'eor\`eme \ref{Alfonso_maintheo1} on utilise alors la densit\'e des classes rationnelles et la semi-continuit\'e de l'ensemble d'Aubry en dimension deux (cons\'equence de \cite{FFR} et \cite{Bernard_Conley}). On obtient que pour toute classe de cohomologie $c$, l'ensemble d'Aubry de $c$ est $\T^2$ tout entier. Un tel ensemble d'Aubry est un graphe lipschtzien au dessus de la section nulle de $T \T^2$, dont la classe de cohomologie est $c$.

En fait on peu dire un peu plus : si un des graphes ci-dessus est feuillet\'e par orbites p\'eriodiques, alors ces orbites p\'eriodiques ont toutes m\^eme p\'eriode. Cela revient \`a exclure les faces radiales de $\beta$. C'est le sens de la   
 \begin{proposition}[\cite{Alfonso}, proposition 3] \label{C1_implies_convex}
Soit  $L$ un lagrangien de  Tonelli autonome sur le tore de dimension deux dont la fonction  $\beta$ est diff\'erentiable en tout  point de $H_1(\T^2,\R)$. Alors quel que soit $h \in H_1(\T^2,\R) \setminus \{ 0 \}$, on a  $R_h = \{ h \}$. En particulier, $\beta_L$ est strictement convexe.
\end{proposition}

\section{Questions}
Le th\'eor\`eme \ref{Alfonso_maintheo1} r\'eduit la question \ref{Burago-Ivanov} \`a celle-ci, pos\'ee dans \cite{Arnaud} :
\begin{question}
Pour un lagrangien de Tonelli, la $C^0$-int\'egrabilit\'e entra\^\i ne-t-elle la $C^1$-int\'egrabilit\'e ?
\end{question}
Un r\'esultat r\'ecent de Sorrentino (\cite{Sorrentino}) affirme que dans le cas d'un hamiltonien de Tonelli, pourvu qu'on puisse trouver suffisement d'int\'egrales pre\-mi\`eres $C^1$, il n'est pas n\'ecessaire de v\'erifier que celles-ci sont en involution.  Il s'agit donc de voir si les graphes lipschitz trouv\'es plus haut sont $C^1$. Pour cela un outil pr\'ecieux est le th\'eor\`eme 2 de \cite{Arnaud}, qui affirme que si le temps un du flot restreint à l'un de ces graphes lipschitz est conjugué, de façon bi-lipschitzienne, à une rotation, alors le graphe en question  est $C^1$. La question \ref{Burago-Ivanov} se ram\`ene donc essentiellement \`a 
\begin{question}
Si la fonction $\beta$ d'un lagrangien de Tonelli sur le tore de dimension deux est $C^1$, l'ensemble de Mather est-il \'egal \`a $\T^2$ tout entier quel que soit la classe de cohomologie ?
\end{question}

\chapter{Conjecture de Ma\~n\'e}
\section{Conjecture forte et conjecture faible}
Dans l'introduction nous affirmons que la conjecture \ref{Mane_forte} contient la conjecture \ref{mane_faible}. 
En effet  \ref{mane_faible} admet une classe  plus large de perturbation, puisqu'en plus d'un potentiel elle autorise \`a rajouter une 1-forme ferm\'ee au lagrangien. Toutefois la pr\'esence d'un ouvert dense dans l'\'enonc\'e emp\^eche \ref{mane_faible} d'\^etre cons\'equence imm\'ediate de \ref{Mane_forte}. Dans \cite{tworemarks} nons d\'emontrons notre affirmation, \`a l'aide des r\'esultats de \cite{CI99} et  \cite{Bernard_Conley}. Pour ce faire on commence par montrer que la conjecture \ref{Mane_forte} est \'equivalente \`a un \'enonc\'e a priori plus fort, o\`u l'ensemble de Mather est remplac\'e par l'ensemble d'Aubry, et o\`u on demande que l'orbite p\'eriodique soit hyperbolique : 
\begin{conjecture}\label{forte_Aubry}
Soient
\begin{itemize}
	\item $M$ une vari\'et\'e ferm\'ee
	\item $L$ un lagrangien de  Tonelli autonome sur  $TM$
	\item $\mathcal{O}(L)$ l'ensemble des $f$ dans  $C^{\infty}(M)$ telles que l'ensemble d'Aubry de $L+f$ consiste en une unique orbite p\'eriodique hyperbolique.
\end{itemize}
Alors l'ensemble  $\mathcal{O}(L)$ est r\'esiduel dans $C^{\infty}(M)$.
\end{conjecture}
Ensuite on d\'emontre que la conjecture \ref{forte_Aubry} contient celle-ci, qui contient clairement  \ref{mane_faible} : 
\begin{conjecture}\label{faible_Aubry}
 Si $L$ est un lagrangien de  Tonelli autonome sur une vari\'et\'e $M$, il existe une partie r\'esiduelle  $\mathcal{O}_4(L)$ de  $C^{\infty}(M)$, telle que pour toute  $f$ dans $\mathcal{O}_2(L)$, il existe un ouvert dense $U(L,f)$ de $H^1 (M,\R)$ tel que, pour tout  $c$ dans $U(L,f)$, l'ensemble d'Aubry de $(L+f,c)$ consiste en une unique orbite p\'eriodique hyperbolique.
\end{conjecture}

La conjecture \ref{Mane_forte} peut \^etre vue comme une version du Closing Lemma adapt\'ee \`a la th\'eorie d'Aubry-Mather. Ceci sugg\`ere qu'elle est vraie en topologie $C^2$ et fausse en topologie $C^k$, pour $k \geq 3$. Rappelons que la topologie $C^2$ sur le hamiltonien correspond \`a la topologie $C^1$ sur le champ de vecteurs. Le cas $C^2$ fait l'objet d'un travail en cours de Figalli et Rifford.  Toutes mes tentatives pour construire un contre-exemple en topologie $C^4$ ont \'echou\'e lamentablement. Je n'ai pu en sauver que l'observation suivante. Supposons que par chance, \'etant donn\'e un lagrangien $L$,  on dispose d'une suite d'orbites p\'eriodiques $\gamma_n$ qui approche l'ensemble de Mather. Alors la premi\`ere id\'ee qui vient \`a l'esprit est de rajouter \`a $\gamma_n$ un potentiel $f_n$ qui s'annule sur  $\gamma_n$, strictement positif partout ailleurs. Si on peut trouver $f_n$ assez grande pour que $\gamma_n$ soit $L+f_n$-minimisante, mais que n\'eanmoins la norme $C^4$ de $f_n$ tende vers z\'ero, on a d\'emontr\'e la partie ''densit\'e'' de la conjecture. Dans \cite{tworemarks} on montre que cette approche naïve ne fonctionne pas. Plus pr\'ecisement, on donne un exemple d'un lagrangien $L$ sur le tore de dimension deux, telle que pour toute orbite p\'eriodique $\gamma$ de $L$, et toute fonction $f$ $C^4$, si $\gamma$ est $L+f$-minimisante, alors la norme $C^4$ de $f$ est born\'ee inf\'erieurement par une constante ne d\'ependant que de $L$.
\section{Transform\'ees de Legendre de classes d'homo\-logie rationnelles}
Puisque la classe d'homologie d'une mesure port\'ee par une orbite p\'eriodique d'un lagrangien p\'eriodique est rationnelle, un premier pas vers la conjecture \ref{faible_Aubry} consiste \`a d\'ecrire l'ensemble des classes de cohomologie qui sont sous-d\'eriv\'ees de $\beta$ en une classe d'homologie rationnelle.

 \begin{theorem}[\cite{ijm}]\label{rationnel_dense_periodique}
 Soient
\begin{itemize}
  \item $M$ une vari\'et\'e ferm\'ee de dimension deux
  \item $L$ un lagrangien de Tonelli autonome sur $TM $.

\end{itemize}
Alors l'ensemble des  classes de cohomologie $c$ qui sont sous-d\'eriv\'ees de $\beta$ en une classe d'homologie rationnelle est  dense dans $H^1(M,\R)$.
\end{theorem} 
Ce    th\'eor\`eme  se d\'eduit, par les th\'eor\`emes 1.1 et 1.2 de \cite{FFR},   de la
\begin{proposition}[\cite{completeproof}]\label{lowdim_dense_1}
Soient 
\begin{itemize}
	\item $M$ une vari\'et\'e ferm\'ee
	\item $L$ un lagrangien de Tonelli sur $TM \times \T$ 
	\item $U$ un ouvert de  $H^1(M,\R)$, tel que pour tout  $c$ dans $U$, $E_c =V_c$ et l'ensemble d'Aubry quotient  $A_c$ est de mesure de Hausdorff unidimensionnelle zero.
\end{itemize}
 Alors la transform\'ee de Legendre $\mathcal{L}(U)$ contient une classe d'homologie rationnelle.
\end{proposition}

\section{R\'eponse \`a la question \ref{mane_faible} en deux degr\'es de libert\'e}
Connaissant  le  th\'eor\`eme  \ref{rationnel_dense_periodique}, on en sait assez pour r\'esoudre le probl\`eme \ref{mane_faible}  dans le cas d'un lagrangien autonome en deux degr\'es de libert\'e. Pour mon embarras \'eternel, dans \cite{ijm} je pr\'etend avoir r\'esolu le probl\`eme \ref{mane_faible} dans ce cas particulier. La preuve est au mieux une esquisse, et elle utilise le r\'esultat de  \cite{AvsM} comme s'il allait de soi. Dans \cite{completeproof} on d\'emontre un \'enonc\'e un peu plus fort  :
 \begin{theorem}\label{Mane}
Soient
\begin{itemize}
  \item $M$ une vari\'et\'e ferm\'ee de dimension deux
  \item $L$ un lagrangien de Tonelli autonome sur  $TM$.
  \end{itemize}
Alors il existe une partie r\'esiduelle  $\mathcal{O}(L)$ de $C^{\infty}(M)$, telle que pour tout $f \in \mathcal{O}(L)$, il existe un ouvert dense $U(L,f)$ de $H^1(M,\R)$, tel que pour tout  $c \in U(L,f)$, l'ensemble d'Aubry  $\mathcal{A}(L+f,c)$ consiste en exactement une orbite p\'eriodique hyperbolique.
\end{theorem} 
Un r\'esultat analogue est d\'emontr\'e pour les lagrangiens p\'eriodiques sur le cercle dans \cite{Osvaldo_09}. 
La diff\'erence avec \cite{ijm} est que l'\'enonc\'e concerne les ensembles d'Aubry plut\^ot que de Mather. L'\'enonc\'e sur les ensembles d'Aubry semble plus utile  au vu des r\'esultats de\cite{A03}, \cite{AIPS05}, et \cite{Bernard07}. 

Voici les grandes lignes de la preuve : 
Par \cite{Mane96}, si $h$ est une classe d'homologie fix\'ee, pour un lagrangien g\'en\'erique il existe une unique mesure $h$-minimisante. L'ensemble des classes d'homologie rationnelles \'etant d\'enom\-brable, pour un lagrangien g\'en\'erique, pour toute classe d'homologie rationnelle $h$, il existe une unique mesure $h$-minimisante. En vertu de la proposition \ref{rational}, cette mesure est port\'ee par une r\'eunion finie d'orbites p\'eriodiques. D'apr\`es \cite{CI99} on peut supposer que ces orbites p\'eriodiques sont hyperboliques. L'ensemble $U$ des classes de cohomologie dont l'ensemble d'Aubry consiste en une orbite p\'eriodique hyperbolique est ouvert \`a cause de l'hyper\-bolicit\'e et de la semi-continuit\'e de l'ensemble d'Aubry en dimension deux. Pour montrer que $U$  est dense, il suffit de montrer qu'il est dense dans les transform\'ees de Legendre de classes rationnelles, en vertu du th\'eor\`eme \ref{rationnel_dense_periodique}. Prenons une classe de cohomologie $c$ qui est sous-d\'eriv\'ee de $\beta$ en une classe rationnelle $h$. L'unique mesure $h$-minimisante est port\'ee par une r\'eunion finie d'orbites p\'eriodiques hyperboliques $\gamma_1, \ldots \gamma_k$. Soit $h_1$ la classe d'homologie de la mesure de probabilit\'e \'equir\'epartie sur $\gamma_1$. Alors $c$ est sous-d\'eriv\'ee de $\beta$ en $h_1$, donc $c$ est approch\'ee par des classes de cohomologie situ\'ees dans l'int\'erieur relatif de $\mathcal{L}(h_1)$. Il reste \`a d\'eterminer l'ensemble d'Aubry $\mathcal{AL}(h_1)$. Si $h_1$ est singuli\`ere, quitte \`a rajouter un potentiel arbitrairement petit pour exclure les orbites homoclines, on montre que $\mathcal{AL}(h_1)$ se r\'eduit \`a un point fixe hyperbolique. Si $h_1$ n'est pas singuli\`ere, on utilise le th\'eor\`eme \ref{AvsM} : $\mathcal{AL}(h_1)$ est \'egal \`a l'ensemble de Mather de $R(h_1)$, la plus grande face radiale de $\beta$ contenant $h_1$. Si $R(h_1)=\left\{h_1 \right\}$, alors $\mathcal{AL}(h_1)$ ne contient que $\gamma_1$, et la preuve est termin\'ee. Si en revanche $h_1$ fait partie d'une face radiale non triviale, il faut utiliser un argument d'approximation  (\cite{AvsM}, lemme C.3). On montre ainsi la densit\'e de $U$, et le th\'eor\`eme.

\section{Application}
Rappelons qu'en vertu du    Th\'eor\`eme 1 de \cite{Bernard07}, lorsque l'ensemble d'Aubry est une r\'eunion finie d'orbites p\'eriodiques hyperboliques, il existe une sous-solution de l'\'equation de Hamilton-Jacobi qui est aussi r\'eguli\`ere que le hamiltonien. Par cons\'equent
\begin{corollary}
Soient
\begin{itemize}
 \item $M$ une vari\'et\'e ferm\'ee de dimension deux
  \item $L$ un lagrangien de Tonelli autonome $C^k$ sur  $TM$, avec $2 \leq k \leq \infty$.
  
  \end{itemize}
  Alors il existe une partie r\'esiduelle  $\mathcal{O}(L)$ de $C^{\infty}(M)$, telle que pour tout $f \in \mathcal{O}(L)$, il existe un ouvert dense $U(L,f)$ de $H^1(M,\R)$, tel que pour tout  $c \in U(L,f)$, il existe une sous-solution $C^k$ de l'\'equation de Hamilton-Jacobi  associ\'ee \`a $(L+f,c)$.
\end{corollary}
\subsection{Question : mesure de l'ouvert dense}
Un r\'esultat r\'ecent de Fathi  (\cite{Fathi_Kyoto}) affirme que,  pour un lagrangien autonome sur une surface, il ne peut y avoir de sous-solution lisse de l'\'equation de Hamilton-Jacobi que si l'ensemble d'Aubry est un tore quasi p\'eriodique ou une r\'eunion d'orbites p\'eriodiques. Un ouvert dense est donc le mieux que l'on puisse esp\'erer en g\'en\'eral dans le corollaire ci-dessus. Il reste \`a d\'eterminer si l'ouvert dense en question peut \^etre de mesure pleine. Le th\'eor\`eme KAM exclut que ce soit le cas  g\'en\'eral. A l'oppos\'e, McShane et Rivin  (\cite{McShane-Rivin}) affirment que c'est le cas lorsque le lagrangien est une m\'etrique de courbure constante n\'egative sur un tore de dimension deux \'epoint\'e. Comme pour le probl\`eme de diff\'erentiabilit\'e, on peut se poser deux questions : 
\begin{itemize}
	\item si le flot d'Euler-Lagrange est d'Anosov, l'ouvert dense est-il de mesure pleine ?
	\item pour un lagrangien g\'en\'erique, l'ouvert dense est-il de mesure pleine ?
\end{itemize}
\chapter{Conjectures de Ma\~n\'e en codimension un}
En th\'eorie de codimension un, le lagrangien est une fonction satisfaisant aux hypoth\`eses suivantes, introduites par Moser dans \cite{Moser}  : 
\begin{itemize}
\item  (H1) : $L\in C^{l,\gamma}(\R^{2n+1})$, $l\ge 2$, $\gamma> 0$.

\item (H2) : $L$ est $1$-p\'eriodique en  $x_1,\dots,x_n,u$.

\item (H3) : il existe  $\delta>0$ tel que 
$$\delta I\le\frac{\partial^2L}{\partial p_i\partial 
p_j}\le{\frac{1}{\delta}}I$$
o\`u  $I$ est la matrice identit\'e de taille  $n$.

\item (H4) : il existe  $C>0$ tel que 
$$\left\vert\frac{\partial^2L}{\partial p\partial x}\right\vert+
\left\vert\frac{\partial^2L}{\partial p\partial u}\right\vert
\le C(1+|p|)$$
$$\left\vert\frac{\partial^2L}{\partial x\partial x}\right\vert+
\left\vert\frac{\partial^2L}{\partial u\partial x}\right\vert+
\left\vert\frac{\partial^2L}{\partial u\partial u}\right\vert
\le C(1+|p|^2) .  $$

\end{itemize}
L'exemple le plus simple est de la forme 
	\[ L(x,u, \nabla u)= \frac{1}{2}\left|\nabla u(x)\right|^2 +f(x,u)
\]
o\`u $f \in C^{l,\gamma}(\R^{n+1})$ est $\Z^{n+1}$-p\'eriodique. 

Les objets que nous recherchons sont des fonctions $u$ de $\R^n$ dans $\R$ qui minimisent localement l'int\'egrale 
\begin{equation}\label{integral_to_be_minimized}
	\int_{\R^n}L(x,u,\nabla(u) )dx
\end{equation}
Moser a d\'emontr\'e l'existence de solutions v\'erifiant une propri\'et\'e suppl\'ementaire : leurs graphes se projettent injectivement dans le tore $\R^{n+1}/\Z^{n+1}$. Les graphes de ces solutions restent \`a distance finie d'un hyperplan non vertical de $\R^n \times \R$ ; Moser a, de plus, d\'emontr\'e que pour tout hyperplan non vertical $\rho$ de $\R^n \times \R$, il existe une solution $u \in  C^{l,\gamma}(\R^{n})$ dont le graphe reste \`a distance finie de 
$\rho$. On dit alors que $\rho$ est la pente moyenne de $u$. Par commodit\'e on identifiera $\rho$ avec l'unique vecteur orthogonal \`a $\rho$, dont la $n+1$-i\`eme coordonn\'ee vaut un. On d\'efinit l'action moyenne d'une telle solution par la formule
$$\lim_{R\to\infty}{\frac{1}{|B(0,R)|}}
\int_{B(0,R)}L(x,u,\nabla u)dx$$
o\`u $B(0,R)$ est la boule euclidienne de rayon $R$, centr\'ee en l'origine, et $|B(0,R)|$ est son volume euclidien. De plus la limite ci-dessus ne d\'epend que de $L$ et  de $\rho$, on la notera donc $\beta_L (\rho)$. La fonction $\rho \longmapsto \beta_L (\rho)$ est surlin\'eaire et strictement convexe (\cite{Senn91}).  

Soient maintenant 
\begin{itemize}
  \item une solution de Moser $u$, de pente moyenne $\rho$
  \item un vecteur   $p=(p_1,\ldots p_n) \in\R^n$,  
  \item une $n$-forme diff\'erentielle lisse $\omega$ sur $\T^{n+1}$.
\end{itemize}
On note $\omega(x,u)\cdot (p,1)$ l'\'evaluation de la  $n$-forme 
$\omega$ sur le  $n$-vecteur (en colonnes)
$$
\begin{array}{cccc}
1 & 0 & \ldots & 0 \\
0 & 1 & \ldots & 0 \\
\ldots  &\ldots & \ldots & \ldots \\
p_1 &  p_2 & \ldots & p_n     
\end{array}
 $$
et il est d\'emontr\'e dans \cite{Bessi09} que la limite suivante existe :     
\begin{equation}\label{corri}
T_u (\omega) := \lim_{R\to\infty}\frac{1}{|B(0,R)|}\int_{B(0,R)}
\omega(x,u(x))\cdot(\nabla u(x),1)dx  , 
\end{equation}
elle d\'efinit alors un $n$-courant $T_u$, de surcro\^\i t 
  $T_u$ est ferm\'e et d'homologie $(\rho,1)$. Les courants ainsi d\'efinis jouent, en th\'eorie de codimension un, le r\^ole des mesures minimisantes de la th\'eorie de Mather. Ils permettent d'\'enoncer les conjectures de Mather en codimension un, sous la forme des questions  \ref{question codim un}.

Dans \cite{codim1} nous r\'epondons par l'affirmative aux questions \ref{question codim un}. Il est piquant de constater que les conjectures de Ma\~n\' e se pr\^etent aussi bien au cadre de codimension un ; une explication possible est que Ma\~n\' e fondait son intuition sur le cas des twists maps, qui \`a certains \' egards est plus proche de la codimension un que de la dimension un. Notons enfin que nous ne disons rien de la conjecture \ref{Mane_forte}, toujours ouverte m\^eme dans le cas des twist maps, du moins en topologie $C^{\infty}$ (voir  \cite{FR}  pour des cas de moindre r\' egularit\' e).

\appendix
\chapter{Preuve du th\'eor\`eme \ref{vmb2_vertices_autonome}}\label{vmb2+}
In \cite{vmb2}, we proved the following
\begin{theorem}
If the $\beta$ function of an autonomous Lagrangian on a manifold $M$ is not differentiable in any direction other than radial 
at  $h$, and the energy level that contains the supports of the $(L,h)$-minimizing measures is supercritical,   
then  $2I_{\R}(h)-1 < \dim M$. 
In particular, if the dimension of $M$ is three, then $I_{\R}(h)\leq 1$. 
\end{theorem}
It turns out that the supercriticality hypothesis is not necessary, as we now prove : 
\begin{theorem}
If the $\beta$ function of an autonomous Lagrangian on a manifold $M$ is not differentiable in any direction other than radial 
at  $h$,    
then  $2I_{\R}(h)-1 < \dim M$. 
In particular, if the dimension of $M$ is three, then $I_{\R}(h)\leq 1$. 
\end{theorem}

Let $h_0 \in H_1 (M,\R)$ be such that the tangent cone to $\beta$ at $h_0$ contains no plane. Let $c \in H^1 (M,\R)$ be such that $\langle c,h_0 \rangle=\alpha (c) + \beta (h_0)$, that is, $c$ defines a supporting hyperplane to $\beta$ at $h_0$. Modifying $L$ by a closed one-form if $c \neq 0$, we may assume $c=0$. Then  $V_0$ has codimension one or zero. We now show that the second case only happens when $h_0=0$. For convenience we also assume that $\alpha(0)=0$.
 \begin{lemma}
 We have $\langle V_0, h_0\rangle=0$.
 Therefore, if  $h_0\neq 0$,
 then $V_0$ has codimension one and $V_0$ is the orthogonal of $h_0$.
 \end{lemma}
 \textbf{Proof.}  Since $0 \in H^{1}(M,\R)$ is a subderivative to $\beta$ at $h_0$, and $L$ is autonomous,  by \cite{Carneiro} the support of any $h_0$-minimizing measure is contained in the energy level $0=\alpha(0)$. Since 
 any $c \in F_0$ is also a subderivative to $\beta$ at $h_0$, the support of any $h$-minimizing measure is contained in the energy level $\alpha(c)$ so we have 
	\[0=\alpha(c )\; \forall c \in F_0.
\]
In other words, the faces of $\alpha$ are contained in the level sets of $\alpha$ for an autonomous Lagrangian.
Besides, since $c$ and $0$ are subderivatives to $\beta$ at $h_0$, we have 
\begin{eqnarray*}
\langle c,h_0 \rangle&=& \alpha (c) + \beta(h_0)=\beta(h_0)	\\
\langle 0,h_0 \rangle&=& \alpha (0) + \beta(h_0)=\beta(h_0)	
\end{eqnarray*}
 whence 
	\[\langle c,h_0 \rangle  =\langle0,h_0\rangle =0 \; \forall c \in F_0.
\]
 $\Box$

 Since $I_{\R}(0)=0$, there is nothing to prove in the case
 $h_0=0$, hence by the lemma above  we may assume $h_0\neq 0$ and $V_0$ has codimension one.

 It would help if we had an ergodic $h_0$-minimizing measure. 
Such a measure need not exist because $(h_0,\beta(h_0))$ may not be an extremal point of the epigraph of $\beta$, 
that is, $h_0$ may  lie in the relative interior of a face of $\beta$. 
Such a face must be radial, i.e. contained in $\left\{ th_0 \co t\in \right[ 0,+\infty \left[ \right\}$, 
for $\beta$ is not differentiable at $h_0$ in any direction but the radial one. 
Then the radial face of $\beta$ containing $h_0$ has at least one non-zero extremal point $th_0$ with $t\in \left] 0,+\infty \right[$. 
Furthermore if the tangent cone to  $\beta$  at $h_0$ is contains no plane, then neither does     the tangent cone to  $\beta$  at $th_0$ (see \cite{AvsM}, Lemma 17).
Since $I_{\R}(th_0)=I_{\R}(h_0)$, we may, for the purpose of proving Theorem 5, assume that $h_0$ itself is an extremal point of $\beta$. 
Then by \cite{Mane92} there exists an ergodic $h_0$-minimizing measure $\mu$.

 Now we proceed with the proof of Theorem 5.
 Let $c_1,\ldots c_{b-1}$ be a basis of $V_0$ and let $\omega_1,\ldots \omega_{b-1}$ be smooth closed one-forms on $M$ 
such that $[\omega_i]=c_i, i=1\ldots k$. 
 Consider the functions
\begin{eqnarray*}
u_{1,i} (x) & :=  &	h(L-\omega_i)\left( x_0,x \right)\\
u_{0} (x) & :=  &	h(L)\left( x_0,x \right).
\end{eqnarray*}	
for $1,\ldots b-1$, and replace each one-form  $\omega_i$ by 
	\[ \omega'_i := \omega_i + du_{1,i} - du_{0}.
\]
This is an almost everywhere defined, integrable and bounded one-form, and it vanishes identically on $\tilde{\mathcal{A}}_0$.

Complete $c_1,\ldots c_{b-1}$ to a basis $c_1,\ldots c_{b}$ of $H^{1}(M,\R)$. Take $\omega_b$ to be any smooth one-form with cohomology $c_b$.
 Pick $a_1,\ldots a_{b}$ an integer basis of $H^{1}(M,\R)$.  Let $(\nu_{ij})_{i,j=1,\ldots b}$ be the matrix 
whose $i$-th line consists of  the coordinates of $a_i$  in the basis $c_1,\ldots c_b$, so 
$a_i = \sum^b _{j=1} \nu_{ij}c_j$. 
Set, for  $j=1, \ldots b$, 
\[
\eta_i := \sum^b_{j=1}\nu_{ij}\omega'_j.
\]
This is an almost everywhere defined, integrable and bounded one-form with cohomology $a_i$.

Now let $(\lambda_{ij})_{i,j=1,\ldots b}$ be the inverse matrix of  $(\nu_{ij})_{i,j=1,\ldots b}$, 
so 
$$
c_i = \sum^b_{j=1}\lambda_{ij}a_j.
$$
 Let $\Lambda$ be the matrix $(\lambda_{ij})_{i=1,\ldots b-1, j=1,\ldots b}$, 
that is, we just delete the last line of $(\lambda_{ij})_{i,j=1,\ldots b}$.
 
 Endow $H_1 (M,\R)$ with the dual basis to $a_1,\ldots a_{b}$ and  $\R^{b-1}$ 
 with its canonical basis. Consider the linear map from $H_1 (M,\R)$ to $\R^{b-1}$ 
whose matrix with respect to said basis is $\Lambda$. We again call this linear map $\Lambda$ for simplicity.
So if $h \in H_1 (M,\R)$, we have 
\begin{eqnarray*}
\Lambda (h) &=& \left(  \sum^b_{j=1}\lambda_{ij} \langle a_j, h \rangle \right)_{i=1,\ldots b-1}\\
&=&  \left(   \langle c_i, h \rangle \right)_{i=1,\ldots b-1}.
\end{eqnarray*}
Thus the kernel of $\Lambda$ is the straight line generated by $h_0$.

 Consider the map 
 \[
\begin{array}{rcl}
\Phi=(\Phi_1,\ldots,\Phi_b) \co M & \longrightarrow & \T^{b}  \\
x & \longmapsto & \left( \int^{x}_{x_0} \eta_i 	\right)_{i=1,\ldots b} \mbox{ mod}\Z^b.
\end{array}
\]
The map $\Phi$ is Lipschitz.
We assume, without loss of generality, that 
$$
\langle a_1,h_0\rangle\neq 0.
$$
 Pick an orbit $\gamma \co  \R \longrightarrow M$ contained in the support $\spt \mu$ of $\mu$.  We have 
	\[ d\Phi (\dot{\gamma}(t))=\left(\eta_i (\dot{\gamma}(t))  \right)_{i=1,\ldots b}
	\mbox{ so } \Lambda \left(d\Phi (\dot{\gamma}(t))  \right)= 
  \left( \omega'_i  (\dot{\gamma}(t))  \right)_{i=1,\ldots b}=0.
\]
 and $\Phi (\gamma (\R))$ is a connected subset of  a straight line
 of direction $h_0$  in $\T^{b}$. 
 Note that $\Lambda \circ \Phi$ is the map
 \[
\begin{array}{rcl}
 \overline{M} &  \longrightarrow  &   \R^{b-1}\\
 x &  \longmapsto &    \left( \int^{x}_{x_0} \omega'_i 	\right)_{i=1,\ldots b} ,
\end{array}
\]
so by Proposition \ref{Holder} $\Lambda \circ \Phi$ is 2-H\"{o}lder on $K_i$. 
 
 We claim that for $\mu$-almost every orbit $\gamma$ the closure of 
 $\Phi (\gamma (\R))$ in $\T^{b}$ is a subtorus of dimension  $I_{\R}(h_0)$.
It is sufficient to prove that the image of $\mu$-almost every orbit is a whole line directed by $h_0$.
This follows from 
\begin{equation}\label{birkhoff}
 \frac{1}{t}\left( \int^{t}_{0} \eta_b ( \dot{\gamma}(t)	)\right) \\
\longrightarrow   \left\langle a_1,h_0 \right\rangle\neq 0
\end{equation}
when $t\rightarrow \pm\infty$,
by  the Birkhoff Ergodic Theorem. Let $U$ be the set
 $$ 
 \big\{ (x,v) \in TM \co \eta_b (x,v) > \frac{1}{2} \left\langle a_1,h_0 \right\rangle \big\}
 $$
 then $U$ is open in $TM$ and by Equation \ref{birkhoff}, we have 
$\mu (U) > \frac{1}{2}$. 

On the other hand $U$  does not meet the zero section of $TM$, so we know from Mather's Graph Theorem
 that the intersection with $U$ of the support of $\mu$ is a Lipschitz lamination whose leaves 
 are the orbits.
In other words, one can cover $\spt \mu \cap U $
with finitely many charts $(U_i,f_i)$ such that:
\begin{itemize}
	\item the $U_i$ are open sets of $M$ that cover $\spt \mu \cap U$
	\item each $f_i$ is a bi-Lipschitz homeomorphism from 
	$]-1,1[^n$ to $U_i$.
		\item The horizontal lines $\{y\}\times ]-1,1[$
		are mapped into pieces of leaf.
		\item The image $\Phi(U_i)$ is contained in a small disk in $\T^b$.
		\end{itemize}
As a consequence, we can lift $\Phi_{U_i}$ to an $\R^b$-valued map.
Let us denote by $K_i$ the set 
of points $y$ in $]-1,1[^{n-1}$ such that 
$f_i(y,0)$ belongs to  the support of $\mu$, so
$f_i^{-1}(\spt \mu)=K_i\times ]-1,1[$.
Let us define a map $\Psi_i:K_i\rightarrow \R^{b-1}$
by the expression
$$
\Psi_i(y)=\Lambda \circ \Phi\circ f_i(y,0).
$$
Since $f_i$ is Lipschitz, $\Psi_i$ is 2-H\"{o}lder on $K_i$.

The push-forward of $\mu$ by $\Phi$ is a measure on a torus of dimension  $I_{\R}(h_0)$, invariant by an irrational flow. Irrational flows on tori are uniquely ergodic so 
$\Phi_* \mu$ is a multiple of the Lebesgue measure. 

Therefore the images $\Phi(\spt \mu \cap U_i)$ 
cover a set of positive Lebesgue measure of a torus of dimension $I_{\R}(h_0)$. Hence
there exists $i$ such that $\Phi(\spt \mu\cap U_i)$ contains a  set of positive Lebesgue measure of a vector space of dimension $I_{\R}(h_0)$. Thus $\Psi_i(K_i)$ contains a  set of positive Lebesgue measure  of a vector space of dimension $I_{\R}(h_0)-1$.
Since $\Psi_i$ is 2-H\"{o}lder on $K_i$, and $K_i$ is  a subset of $\R^{\dim M -1}$, Ferry's Lemma  implies $1/2(\dim M -1)> I_{\R}(h_0)-1$, that is, $\dim M > 2I_{\R}(h_0)-1$.
$\Box$

\chapter{Adh\'erence de courbes sur les tores}\label{reponse_courbe}
Dans cette annexe nous donnons les maigres \'el\'ements dont nous disposons concernant  la question \ref{question_courbe}.
Soient 
\begin{itemize}
	\item $P$ un plan totalement irrationnel de $\R^3$ (i.e. $P \cap \Z^3 = \left\{0\right\}$)
	\item $\gamma$ une courbe continue \`a valeurs dans $P$ (sur laquelle les hypoth\`eses seront pr\'ecis\'ees plus tard ; en g\'en\'eral on aimerait que la courbe ait une classe d'homologie, ou direction,  asymptotique, i.e. que $t\inv \gamma(t)$ ait une limite non nulle). 
\end{itemize}
 On se pose la question suivante : \textit{quelle est l'adh\'erence de $\gamma$ dans $\T^3 = \R^3 / \Z^3$} ? On cherche une condition suffisante pour que ladite adh\'erence soit de dimension de Hausdorff au moins deux.
 Soient 
\begin{itemize}
	\item  $l \in \left(\R^3\right)^{\ast}$ telle que $P= \ker l$
	\item $D_0 = \left[0,1\right]^3$ le cube unit\'e de $\R^3$
	\item $ A:= \overline{\gamma (\R) +\Z^3}$
	\item $ A_0:= A\cap D_0$
	\item $\mbox{Leb}$ la mesure de Lebesgue sur $\R$.
\end{itemize}
Remarquons que 
\begin{itemize}
	\item $A$ est une r\'eunion d\'enombrable de translat\'es de $A_0$
	\item $A_0$ est compact
	\item $\mbox{Leb}\left(l(A_0)\right)>1$ implique $\dim A_0 \geq 2$.
\end{itemize}
La derni\`ere assertion est cons\'equence du fait que $A_0$ est localement hom\'eo\-morphe au produit cart\'esien de $l(A_0)$ par un intervalle.
Soit pour tout $n \in \Z$ $Z_n = (x_n,y_n,z_n) \in \Z^3$ tel que $(D_0 - Z_n)_{n \in \Z}$ est la suite des domaines fondamentaux travers\'es par $\gamma (t)$ pour $t \in \R$.
Le lemme qui suit nous ram\`ene \`a un probl\`eme de dynamique symbolique, lequel, \`a d\'efaut d'\^etre plus facile, a au moins l'avantage de se pr\^eter \`a des simulations en Maple (dont l'auteur tient des images \`a la disposition du lecteur curieux).
\begin{lemma}\label{annexeB}
On a 
	\[l(A_0) = \overline{l(Z_n)_{n \in \Z}}.
\]
\end{lemma}
\textbf{Preuve} : pour tout $n \in \Z$ il existe $t_n \in \R$ tel que $\gamma (t_n) \in D_0 - Z_n$. Alors $\gamma (t_n) + Z_n \in D_0 \cap A = A_0$, donc $l\left(\gamma (t_n) + Z_n \right) = l\left(Z_n \right) \in l(A_0)$. Puisque $A_0$ est compact on en d\'eduit
	\[l(A_0) \supset \overline{l(Z_n)_{n \in \Z}}. 
\]
R\'eciproquement, soit $X \in \stackrel{\circ}{D}_0 \cap A$. Par d\'efinition de $A$ il existe une suite $\gamma (t_k) + M_k$ de $\gamma (\R) +\Z^3$ telle que 
	\[ \lim_{k \rightarrow \infty } \left(\gamma (t_k) + M_k \right) = X.
\]
Puisque $X$ se trouve dans l'int\'erieur de $D_0$, pour $k$ assez grand on a $\gamma (t_k) + M_k \in D_0$, en particulier il existe $n(k)$ tel que $M_k = Z_{n(k)}$. On a alors $l\left(\gamma (t_k) + M_k \right)= l(Z_{n(k)}$ donc $l(X)= \lim l\left(\gamma (t_k) + M_k \right)$, d'o\`u
	\[l\left(\stackrel{\circ}{D}_0 \cap A\right) \subset \overline{l(Z_n)_{n \in \Z}} 
\]
puis 
\[l(A_0) \subset \overline{l(Z_n)_{n \in \Z}}. 
\]
ce qui conclut la preuve du lemme.
$\Box$

Sans perte de g\'en\'eralit\'e on peut supposer que $l(x,y,z)= z - \alpha x -\beta y$ avec $1,\alpha,\beta$ rationnellement ind\'ependants.
Notons que
	\[ 
\mbox{Leb} \left(\overline{l(Z_n)}_{n \in \Z}\right) >0 \Longleftrightarrow 
\mbox{Leb} \left(\overline{\left(\alpha x_n + \beta y_n \mbox{ mod } 1\right)}_{n \in \Z}\right) >0 .
\]
 et de plus, \'etant donn\'ee une suite $(x_n,y_n)$ dans $\Z^2$, telle que $(x_{n+1}, y_{n+1})= (x_n,y_n) \pm (1,0) \mbox{ ou } (0,1)$, il existe  une suite enti\`ere $z_n$ telle que $(x_n,y_n, z_n)$ est la suite des domaines fondamentaux travers\'es par une courbe continue contenue dans $P$.

D\'esormais nous nous int\'eresserons donc \`a la r\'epartition modulo un de  suites $\alpha x_n + \beta y_n$, avec $(x_{n+1}, y_{n+1})= (x_n,y_n) \pm (1,0) \mbox{ ou } (0,1)$. 
Observons qu'une telle suite n'est pas forc\'ement dense modulo un. En effet, supposons que $0 < \alpha < \beta < 0.1$ et soit 
$ A := \min \left\{\alpha,\beta - \alpha \right\}$. Prenons $(x_0,y_0) = (0,0)$ et $(x_{n+1}, y_{n+1})=(x_n, y_n)+ (0,1)$ jusqu'au moment o\`u on arrive dans $\left]-\beta,0 \right[$. Alors, si on est tomb\'e dans $\left]-\beta,-\alpha \right[$, on prend $(x_{n+1}, y_{n+1})=(x_n, y_n)+ (1,0)$, qui tombe dans $\left]-\beta + \alpha, 0 \right[$. Ensuite on prend $(x_{n+2}, y_{n+2})=(x_{n+1}, y_{n+1})+(0,1)$ et on arrive dans $\left]\alpha,\beta\right[$. Alors on joue  $(x_{n+3}, y_{n+3})=(x_{n+2}, y_{n+2})+(0,1)$  et enfin $(x_{n+4}, y_{n+4})=(x_{n+3}, y_{n+3})-(1,0)$ pour remettre \`a z\'ero le compteur des $x_n$.

Si on contraire on est tomb\'e dans $\left]-\alpha,0 \right[$, on joue $(x_{n+1}, y_{n+1})=(x_n, y_n)+ (0,1)$, qui arrive dans $\left]\beta - \alpha, \beta \right[$. 
La suite $\alpha x_n + \beta y_n$ modulo un ainsi construite \'evite l'intervalle $\left]0,A \right[$ et v\'erifie $x_n= 0 \mbox{ ou } 1$, en particulier $x_n / y_n $ tend vers z\'ero.
\section{Suites al\'eatoires}
La difficult\'e \`a construire des suites dont l'adh\'erence serait de dimension de Hausdorff $<1$ s'explique peut-\^etre par le fait qu'un tel ph\'enom\`ene, si jamais il se produit, est rare, comme le montre le lemme suivant.
Dans ce paragraphe on s'int\'eresse \`a des suites al\'eatoires au sens de la mesure de probabilit\'e invariante par le d\'ecalage, et par l'\'echange de $\alpha$ et $\beta$. Notons que par la loi des grands nombres une suite g\'en\'erique compte statistiquement autant de 
$\alpha$ que de  $\beta$. En modifiant la mesure on peut obtenir toutes les fr\'equences possibles d'apparition de  $\alpha$, c'est \`a dire toutes les classes d'homologie possibles pour la courbe $\gamma$ de la question originelle.
\begin{lemma}
Presque toute suite est dense dans le cercle.
\end{lemma}
Il suffit pour d\'emontrer le lemme, de montrer qu'une suite al\'eatoire contient des r\'ep\'etitions arbitrairement longues de $\alpha$, puisque la suite $(n\alpha)_{n  \in \N}$ est dense dans le cercle. 

Soit $n  \in \N$. L'ensemble $E_n$ des suites ayant $n$ $\alpha$ cons\'ecutifs est invariant par le d\'ecalage. Celui-ci \'etant ergodique, $E_n$ est de mesure nulle ou pleine. Mais il ne peut \^etre de mesure nulle, puisque il contient l'ensemble des suites qui commencent par $n$ $\alpha$ cons\'ecutifs, et ce dernier est de mesure $2^{-n}>0$.
$\Box$

On peut se poser une question plus fine : \textit{une suite al\'eatoire est-elle \'equir\'epartie ?}
Il para\^\i t utile, au premier abord, de disposer d'une orbite g\'en\'erique explicite du d\'ecalage. On propose le candidat suivant : on num\'erote les suites  finies de $\alpha$  et $\beta$ (les mots finis sur l'alphabet \`a deux lettres), et on les met bout \`a bout, dans l'ordre croissant de leurs num\'eros. L'orbite de cette suite par le d\'ecalage est clairement dense, mais est-elle  \'equir\'epartie ?
Quelles sont la ou les mesures de comptage d\'efinies sur le cercle par cette suite ? Y en a-t-il d'autres que celle de Lebesgue ?

Il est clair qu'une suite qui contient des r\'ep\'etitions arbitrairement longues d'un m\^eme motif est dense dans le cercle. Pour \'eviter cet inconv\'enient, il m'a \'et\'e sugg\'er\'e par B. Rittaud de consid\'erer une suite invariante par substitution, par exemple la substitution de Fibonacci ($0 \mapsto 1, 1 \mapsto 01$). On est donc amen\'e \`a adopter un autre point de vue sur la question : on fixe maintenant la suite de $0$ et de $1$, et on fait varier les nombres $\alpha$ et $\beta$. Il serait int\'eressant  de relier les propri\'et\'es de la suite \`a des propri\'et\'es arithm\'etiques (par exemple, diophantines) de $\alpha$ et $\beta$. Exp\'erimentalement on constate que pour $(\alpha,\beta)=$
\[
(\sqrt{2},  \sqrt{3}), (\pi, \pi^2), (\pi, \sqrt{2}), (\pi, \pi^3), (\pi^3, \pi), (\pi, e), ( e, \pi),  \left( \sqrt{2}, \frac{1+\sqrt{5}}{2}\right),
\] 
la suite est dense dans le cercle ; en revanche pour 
\[(\alpha,\beta)=\left(  2 \frac{2-\sqrt{3}}{1+\sqrt{5}} , \sqrt{3}  \right), \left( 2 \frac{2-\sqrt{2}}{1+\sqrt{5}}, \sqrt{2} \right)
\]
elle \'evite un intervalle. L'auteur ne conna\^\i t \`a ce jour pas de d\'emonstration rigoureuse de ces faits.

\section{Cas o\`u la courbe est \`a distance born\'ee d'une droite}
Supposons donn\'es
\begin{itemize}
	\item un vecteur $v$ totalement irrationnel de $\R^n$
	\item une courbe continue $\gamma \co \R \longrightarrow \R^n$ telle que $\forall t \in \R,\  \gamma(t).v = 0$
	\item une droite $D$ de $\R^n$ et $R > 0$ tels que $\forall t \in \R,\  \exists x \in D, \|x-\gamma(t)\| \leq R$, et 
	$\forall x \in D,\  \exists t \in \R, \|x-\gamma(t)\| \leq R$.
\end{itemize}
Notons $A$ l'adh\'erence de $\gamma (\R) + \Z^n$ dans $\R^n$.
\begin{proposition}
La dimension de Hausdorff de $A$ est au moins deux.
\end{proposition}
{\bf Preuve.}
Notons qu'on peut supposer $\forall x \in D,\ x.v =0$. En effet, soient $x \in D$, et $t \in \R$ tels que $\|x-\gamma(t)\| \leq R$. Alors 
$|x.v| = |(x-\gamma(t)).v| \leq \|v\|R$, dons $x \mapsto x.v$ est born\'e sur $D$. Comme $D$ est une droite, $x \mapsto x.v$ est constant sur $D$. Quitte \`a changer $R$, on peut translater $D$ pour que  $x \mapsto x.v$ soit identiquement nulle sur $D$. 
Notons
\begin{itemize}
	\item $V_C$ le voisinage tubulaire de rayon $C$ de $D$, de sorte que $\gamma(\R) \subset V_R$
	\item $D_0 (C)$ la boule  ferm\'ee de centre l'origine et de rayon $C$
	\item $E := \{ z \in \Z^n \co D +z \subset V_R \}$.
\end{itemize}
Alors $\forall z \in E$, $\forall t \in \R$, $\gamma (t) + z \in V_{2R}$. Puisque $v$ est totalement irrationnel, $E.v$ est dense dans un intervalle $I$. Soit $\theta \in I$. Il existe une suite $z_n \in E$ telle que $v.z_n \longrightarrow \theta$. Puisque $z_n \in E$ pour tout $n$, il existe une suite $x_n \in D$ telle que $x_n+z_n \in B_0(R)$. Donc il existe une suite $t_n \in \R$ telle que $\gamma(t_n) +z_n \in B_0(2R)$, donc $v.z_n \in v.\left(A \cap B_0(2R)\right)$. Puisque $A \cap B_0(2R)$ est compact, cela implique $\theta \in v.\left(A \cap B_0(2R)\right)$, donc 
$I \subset v.A$. Comme $x \mapsto x.v$ est Lipschitz, et constante sur chaque translat\'e de $\gamma$, cela entra\^\i ne que $A$ est de dimension de Hausdorff  au moins deux.
$\Box$
\section{Quasi-cristaux}
Pour les besoins de la cause on appelera quasi-cristal associ\'e \`a $P$, l'ensemble suivant :
	\[ \mathcal{QC}(P) := \left\{ (x,y,z) \in \Z^3 \co z < \alpha x + \beta y < z+1 \right\}. 
\]
On dit que deux points $(x,y,z)$ et $(x',y',z')$ de $\mathcal{QC}(P)$ sont voisins si $\left|x-x'\right|+\left|y-y'\right|+\left|z-z'\right|=1$. On appelle hauteur d'un point $(x,y,z)$ de $\mathcal{QC}(P)$ la diff\'erence $\alpha x + \beta y -z$.  On appelle cha\^\i ne dans $\mathcal{QC}(P)$ une suite $(x_n,y_n,z_n)$ dans $\mathcal{QC}(P)$, telle que $(x_n,y_n,z_n)$ est voisin de $(x_{n+1},y_{n+1},z_{n+1})$ pour tout $n$. On dit qu'une cha\^\i ne infinie $(x_n,y_n,z_n)$ admet une direction asymptotique  si $n\inv (x_n,y_n,z_n)$ admet une limite non nulle. On dit qu'une partie de $\mathcal{QC}(P)$ est connexe si deux points quelconques sont joints par une cha\^\i ne finie. Le lemme \ref{annexeB} sugg\`ere donc la  question  : 
\begin{question}
Etant donn\'ee une cha\^\i ne infinie $(x_n,y_n,z_n)$, admettant une direction asymptotique, 
 la suite des hauteurs  de $(x_n,y_n,z_n)$ est-elle de mesure positive ?  dense dans un intervalle ?
\end{question}
On peut \^etre tent\'e  de prendre le probl\`eme \`a l'envers : si $K$ est une partie de $\left[0,1\right]$, on appelle $\mathcal{QC}(P,K)$ l'ensemble des points de $\mathcal{QC}(P)$ dont la hauteur n'est pas dans $K$. On se demande alors s'il existe une taille critique pour $K$ au del\`a de laquelle les composantes connexes de $\mathcal{QC}(P,K)$ ne peuvent contenir de cha\^\i ne infinie admettant une direction asymptotique. Naturellement on peut prendre pour $K$ un intervalle suffisement grand ; la question devient int\'eressante si on suppose que les composantes connexes de $K$ sont petites (par exemple, plus petites que $\left|\alpha-\beta\right|$). En fait, au vu de simulations, on a m\^eme envie de demander s'il existe une taille critique pour $K$ au del\`a de laquelle les composantes connexes de $\mathcal{QC}(P,K)$ sont finies.  Dans les  figures qui suivent, on prend d'abord pour $K$ la plus grande composante connexe du compl\'ementaire de l'ensemble de Cantor triadique, puis les trois plus grandes, puis les sept, quinze, et trente-et-une plus grandes. 

\begin{figure}[h]
\hspace{1.5cm}
\includegraphics[scale = 0.5]{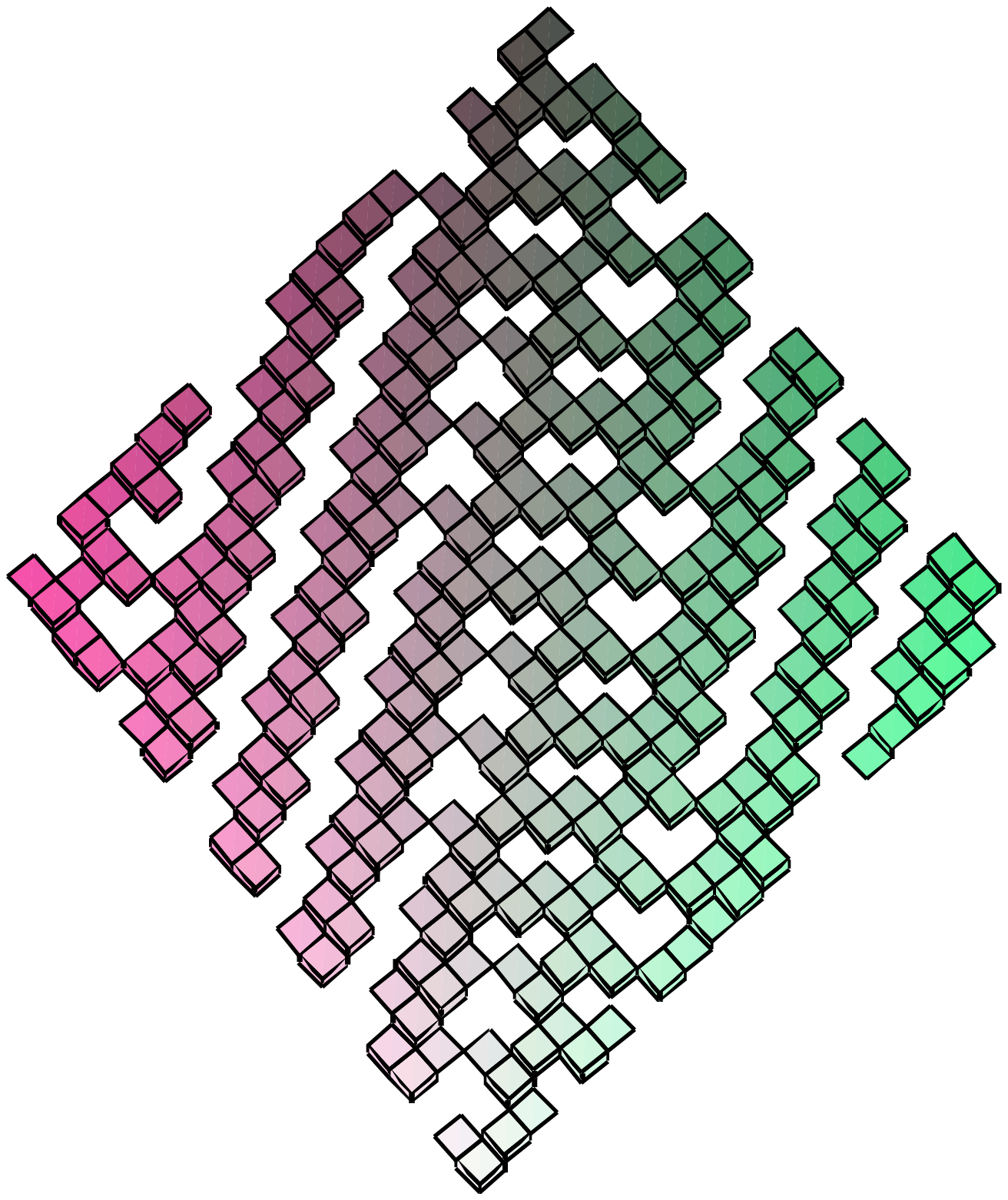} 
\end{figure}
\begin{figure}[h]
\includegraphics[scale = 0.5]{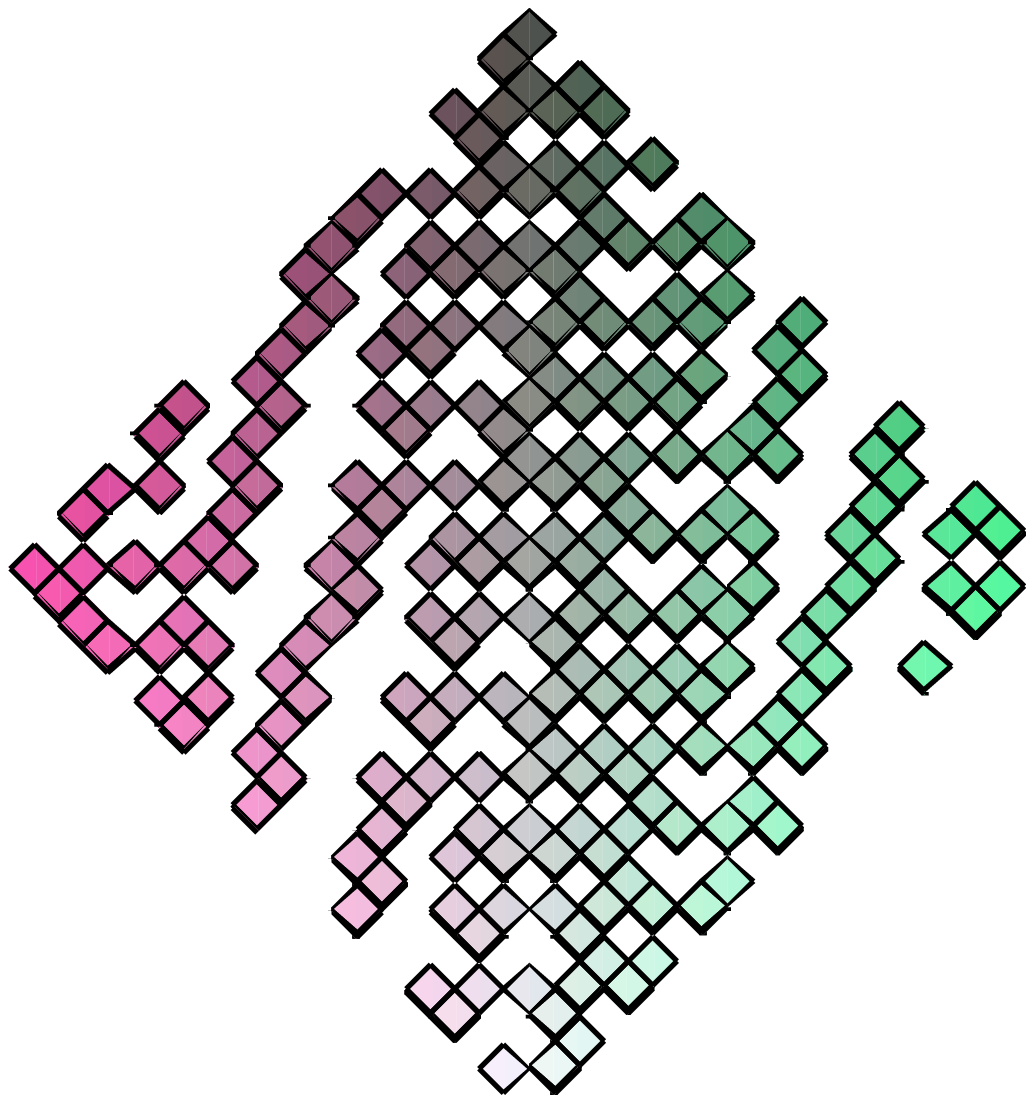} 
\includegraphics[scale = 0.5, angle = -13]{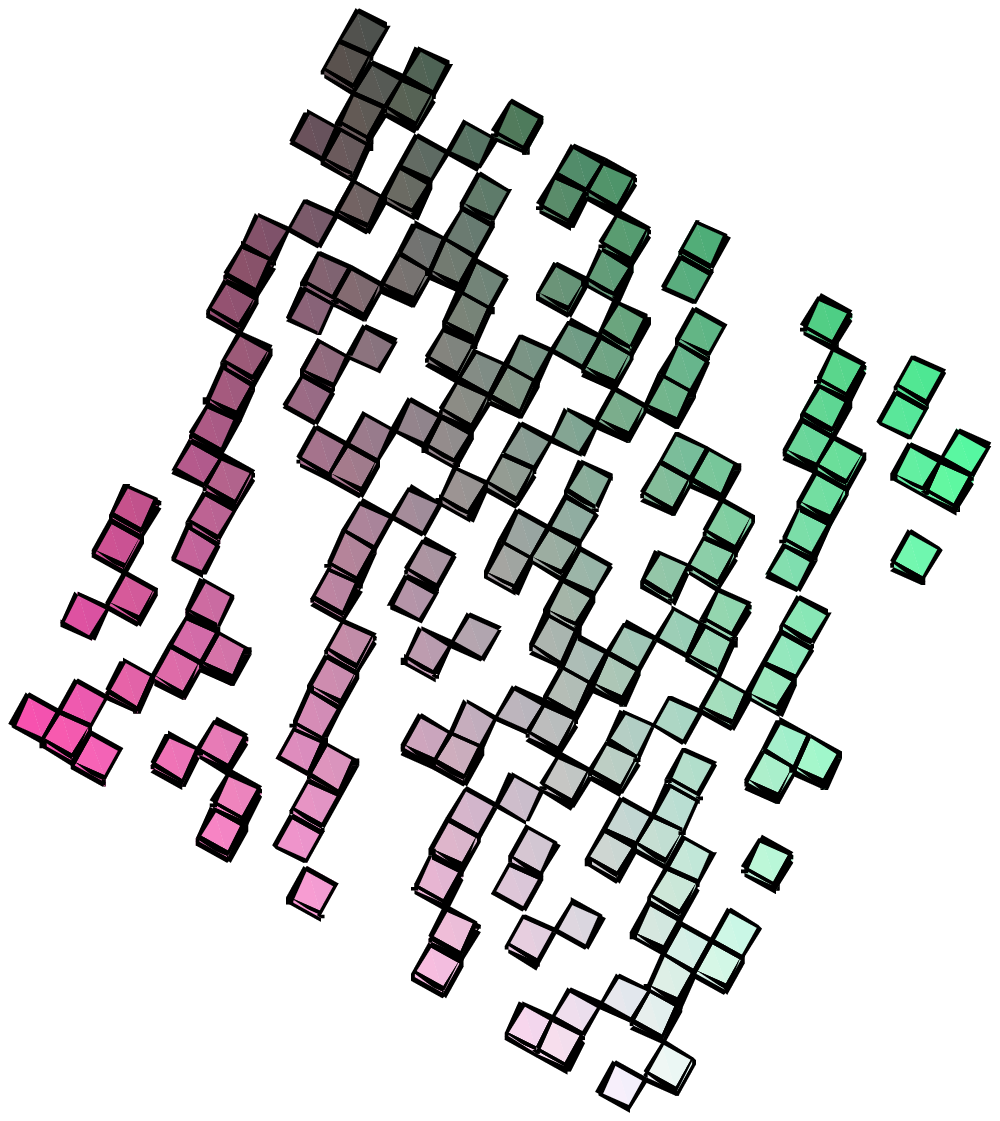} 
\includegraphics[scale = 0.5]{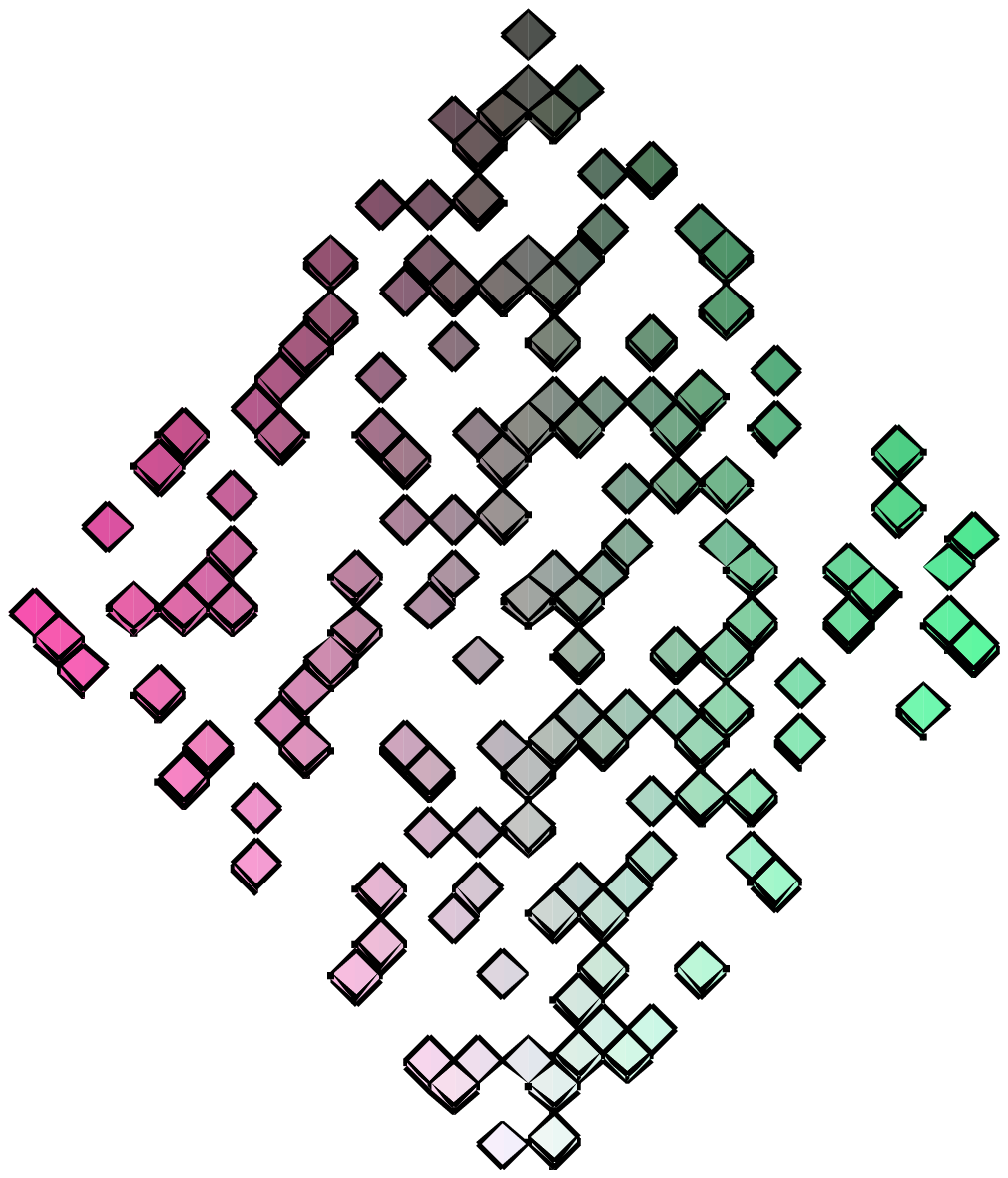} 
\includegraphics[scale = 0.5, angle = 10]{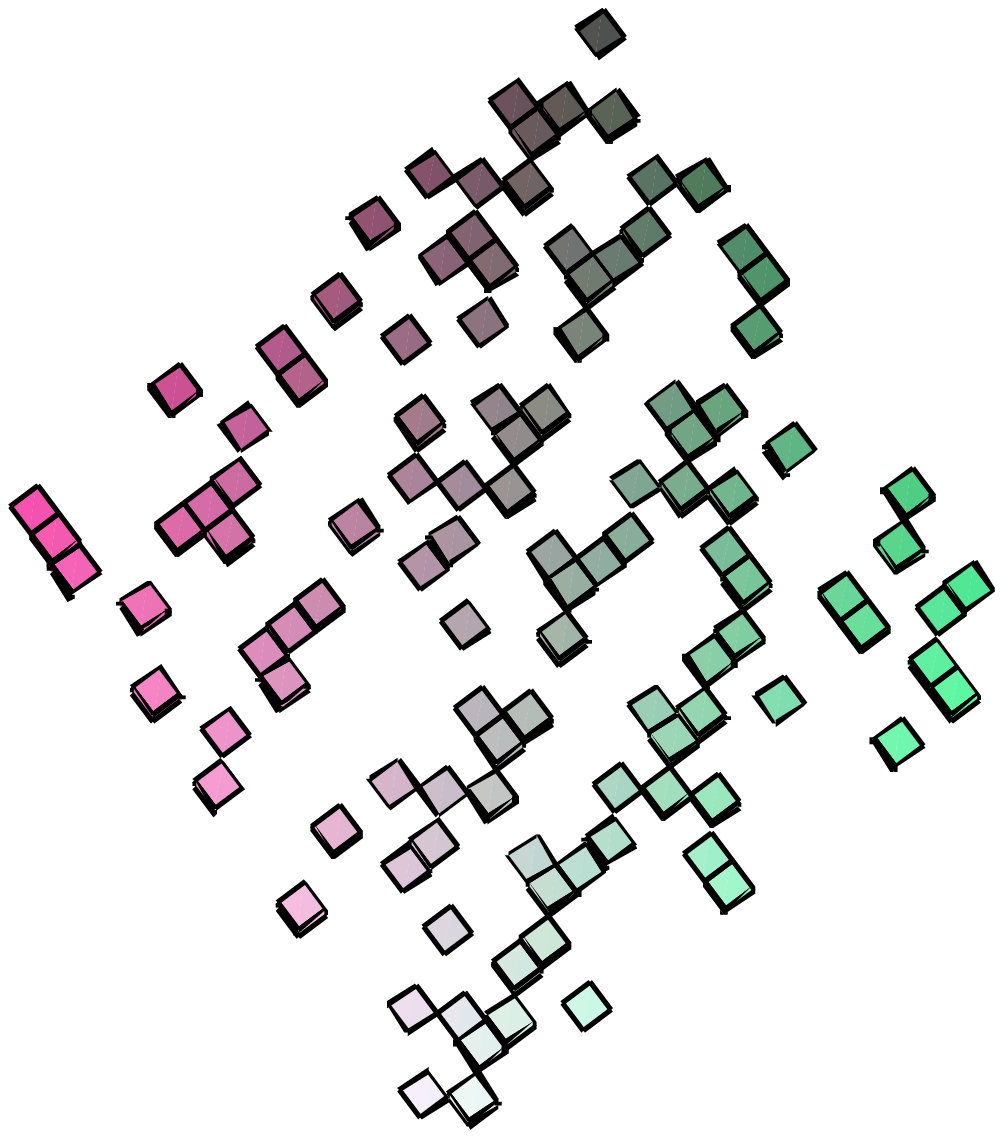} 
\end{figure}
%


{\small
\begin{thebibliography}{99}
\bibitem[Mt96]{these} D. Massart
{\em Normes stables des surfaces} \\
thèse de doctorat, Ecole Normale Supérieure de Lyon, 1996
%
\bibitem[Mt97a]{cras} D. Massart
{\em Normes stables des surfaces} \\
C.R.A.S. S\'erie 1 \textbf{324} (1997), no. 2, 221-224
%
\bibitem[Mt97]{gafa} D. Massart
{\em Stable norms of surfaces: local structure of the unit ball at rational directions}  \\
Geom. Funct. Anal.  \textbf{7}  (1997),  no. 6, 996--1010.
%
\bibitem[Mt03]{ijm} D. Massart 
{\em On Aubry sets and Mather's action functional}\\ 
 Isra\"el Journal of Mathematics {\bf 134} (2003), 157-171.
%
\bibitem[Mt07]{soussol} D. Massart 
{\em Subsolutions of time-periodic Hamilton-Jacobi equations}\\
 Ergodic Theory and Dynamical Systems   \textbf{27} (2007), no. 4, 1253-1265.
%
\bibitem[BaM08]{nonor} F. Balacheff, D. Massart {\em Stable norms of non-orientable surfaces}\\
  Ann. Inst. Fourier (Grenoble)  \textbf{58}  (2008),  no. 4, 1337--1369.
%
\bibitem[Mt09]{vmb2} D. Massart 
{\em Vertices of Mather's beta function, II}\\
 Ergodic Theory and Dynamical Systems \textbf{29} (2009), no. 4,  1289-1307 
%
\bibitem[Mt10]{tworemarks} D. Massart 
{\em Two remarks about Ma\~n\'e's conjecture}\\
Regular and Chaotic Dynamics, 2010, \textbf{15}, No. 6, pp. 646-651.
%
\bibitem[Mta]{AvsM} D. Massart 
{\em Aubry sets vs Mather sets in two degrees of freedom}\\ 
preprint arXiv:0803.2647 [math.DS]\\
 \`a para\^\i tre dans  Calculus of Variations and Partial Differential Equations
%
\bibitem[MS]{Alfonso} D. Massart, A. Sorrentino 
{\em Differentiability of Mather's average action and integrability on closed surfaces }\\
preprint 	arXiv:0907.2055 [math.DS]
%
\bibitem[BeM]{codim1} U. Bessi, D. Massart 
{\em Ma\~n\'e's conjectures in codimension one}\\
 preprint 	arXiv:1009.5474v1 [math.AP]\\ 
%
\bibitem[Mtb]{completeproof} D. Massart 
{\em Generic Aubry sets in two degrees of freedom}\\
travail en cours
%
\newpage
%

\textbf{Articles cit\'es en r\'ef\'erence}
\vspace{0.5cm}
\bibitem[A03]{A03}
N. Anantharaman 
{\em Counting geodesics which are optimal in homology}\\
Ergodic Theory Dynam. Systems 23 (2003), no. 2, 353--388. 
%
\bibitem[AIPS05]{AIPS05}
N. Anantharaman ; R. Iturriaga ; P. Padilla ; H. S\'anchez-Morgado \\
{\em Physical solutions of the Hamilton-Jacobi equation}\\
Discrete Contin. Dyn. Syst. Ser. B 5 (2005), no. 3, 513--528. 
%
\bibitem[Ar08]{Arnaud}
M. C.  Arnaud.
{ \em Fibr\'es de Green et r\'egularit\'e des graphes $C^0$-Lagrangiens invariants par un flot de Tonelli}
 Ann. Henri Poincar\'e, 9 (5): 881--926, 2008.
%
\bibitem[A64]{Arnold}
V. I. Arnol'd 
{\em Instability of dynamical systems with many degrees of freedom}
Dokl. Akad. Nauk SSSR 156 1964 9--12. 
%
\bibitem[AB06]{Bangert_Auer} F. Auer, V. Bangert {\em Differentiability of the stable norm in codimension one}
  Amer. J. Math.  128  (2006),  no. 1, 215--238.
%
\bibitem[BB06]{BB06}I. K. Babenko, F. Balacheff {\em Sur la forme de la boule unit\'e de la norme stable unidimensionnelle}   Manuscripta Math.  119  (2006),  no. 3, 347--358.
%
\bibitem[Ba87]{Bangert unique}V. Bangert {\em A uniqueness theorem for $\Z^n$-periodic variational problems}, Comment. Math. Helv., 62 (1987), 511-531.
%
 \bibitem[Ba88]{Bangert_twist} V. Bangert, {\em Mather sets for twist maps and geodesics on tori}  Dynamics reported, Vol. 1,  1--56, Dynam. Report. Ser. Dynam. Systems Appl., 1, Wiley, Chichester, 1988.
 %
\bibitem[Ba89]{Bangert} V. Bangert {\em On minimal laminations of the torus}, Ann. Inst. Henri Poincar\'e\  6 (1989), no 2, 95--138.
%
\bibitem[Ba90]{Bangert90} V. Bangert {\em Minimal geodesics}  Ergodic Theory Dynam. Systems  10  (1990),  no. 2, 263--286.
%
\bibitem[Ba94]{Bangert94} V. Bangert {\em Geodesic rays, Busemann functions and monotone twist maps}  Calc. Var. Partial Differential Equations  2  (1994),  no. 1, 49--63.
%
\bibitem[Be02]{Bernard_Fourier} P. Bernard {\em Connecting orbits of time dependent Lagrangian systems}  Ann. Inst. Fourier (Grenoble)  52  (2002),  no. 5, 1533--1568.
%
%
\bibitem[Be07]{Bernard07}  P. Bernard {\em  	 
Smooth critical sub-solutions of the Hamilton-Jacobi equation }  \\
Math. Res. Lett.  14  (2007),  no. 3, 503--511. 
%
\bibitem[Be08]{be08} P. Bernard {\em The dynamics of pseudographs in convex Hamiltonian systems} J. Amer. Math. Soc. 21 (2008), no. 3, 615--669. 
%
\bibitem[Be]{Bernard_Conley} 	Patrick Bernard
{\em  On the Conley Decomposition of Mather sets } to appear, Revista Iberoamericana de Matem\'aticas
%
\bibitem[BB07]{Bernard-Buffoni}
P. Bernard, B. Buffoni {\em Optimal mass transportation and Mather theory} J. Eur. Math. Soc. 9 (2007), 85-121.
%
\bibitem[BC08]{Bernard-Contreras}
 P. Bernard, G.  Contreras {\em  A generic property of families of Lagrangian systems}  Ann. of Math. (2)  167  (2008),  no. 3, 1099--1108.
 %
%
\bibitem[Be09]{Bessi09} U. Bessi 
{\em 
 Aubry sets and the differentiability of the minimal average action in codimension one} ESAIM Control Optim. Calc. Var. \textbf{15} (2009), no. 1, 1--48. 
%
\bibitem[BI94]{Burago-Ivanov} D. Burago, S. Ivanov {\em Riemannian tori without conjugate points are flat}  Geom. Funct. Anal.  4  (1994),  no. 3, 259--269.
%
\bibitem[BIK97]{BIK} D. Burago, S. Ivanov, B. Kleiner {\em On the structure of the stable norm of periodic metrics}  Math. Res. Lett.  4  (1997),  no. 6, 791--808.
%
\bibitem[C95]{Carneiro} M. J. D. Carneiro. {\em On minimizing measures of the action of autonomous Lagrangians}  Nonlinearity 89 (1995), 1077--1085.
%
\bibitem[CL99]{CL99}
M.J. Carneiro, A.  Lopes {\em On the minimal action function of autonomous Lagrangians associated to magnetic fields}  Ann. Inst. H. Poincar\'e Anal. Non Lin\'eaire  16  (1999),  no. 6, 667--690.
%
\bibitem[CI99]{CI99} G. Contreras, R. Iturriaga {\em Convex Hamiltonians without conjugate points}  \\
Ergodic Theory Dynam. Systems  19  (1999),  no. 4, 901--952.
%
%
\bibitem[CMP04]{CMP} G. Contreras, L. Macarini, G. Paternain, {\em Periodic orbits for exact magnetic flows on surfaces}    Int. Math. Res. Not.  2004,  no. 8, 361--387. 
%
%
\bibitem[F]{Fathi_bouquin} A. Fathi 
{\em Weak KAM theorem in Lagrangian dynamics} to appear, Cambridge University Press.
%
\bibitem[F09]{Fathi_Kyoto}
A. Fathi.
 {\em Denjoy-Schwartz and Hamilton-Jacobi}
 RIMS meeting: Viscosity solutions of differential equations and related topics, Kyoto University, 2008.

%
\bibitem[FFR09]{FFR} A. Fathi, A. Figalli, L. Rifford
{\em On the Hausdorff Dimension of the Mather Quotient}  Comm. Pure Appl. Math.  62  (2009),  no. 4, 445--500. 
%
\bibitem[FS04]{FS} A. Fathi, A.  Siconolfi {\em Existence of $C\sp 1$ critical subsolutions of the Hamilton-Jacobi equation}  Invent. Math.  155  (2004),  no. 2, 363--388.
%
\bibitem[F75]{Ferry} S.Ferry, {\em  When $\epsilon $-boundaries are manifolds}  Fund. Math.  90  (1975/76), no. 3, 199--210.
%
\bibitem[FR]{FR} A. Figalli, L. Rifford, {\em Closing Aubry sets} preprint 
%
\bibitem[GLP81]{Gromov-Lafontaine-Pansu}
M. Gromov
{\em Structures m\'etriques pour les vari\'et\'es riemanniennes}
Edited by J. Lafontaine and P. Pansu. Textes Math\'ematiques [Mathematical Texts], 1. CEDIC, Paris, 1981. 
%
\bibitem[H39]{Hedlund} G. Hedlund, {\em   Geodesics on a two-dimensional Riemannian manifold with periodic coefficients}  Ann. of Math. (2)  33  (1932),  no. 4, 719--739
%
%
\bibitem[Mn92]{Mane92} 
R. Ma\~n\'e {\em On the minimizing measures of Lagrangian dynamical systems}  Nonlinearity  5  (1992),  no. 3, 623--638.
%
\bibitem[Mn95]{Mane95} 
R. Ma\~n\'e
{\em Ergodic variational methods: new techniques and new problems}  Proceedings of the International Congress of Mathematicians, Vol. 1, 2 (Z\H{u}rich, 1994),  1216--1220, Birkh\H{a}user, Basel, 1995.
%
\bibitem[Mn96]{Mane96} 
R. Ma\~n\'e
{\em Generic properties and problems of
minimizing measures of Lagrangian systems} 
Nonlinearity {\bf 9} (1996), no. 2, 273--310.
%
\bibitem[Mn97]{Mane97} 
R. Ma\~n\'e
{\em Lagrangian flows: the dynamics of globally minimizing orbits}  Bol. Soc. Brasil. Mat. (N.S.)  28  (1997),  no. 2, 141--153.
%
%
\bibitem[Mr90]{Mather90} J. N. Mather {\em Differentiability of the minimal average action as a function of the rotation number}  Bol. Soc. Brasil. Mat. (N.S.)  21  (1990),  no. 1, 59--70. 
%
\bibitem[MrF91]{Mather-Forni} J. N. Mather, G. Forni {\em Action minimizing orbits in Hamiltonian systems}  Transition to chaos in classical and quantum mechanics (Montecatini Terme, 1991),  92--186, Lecture Notes in Math., 1589, Springer, Berlin, 1994. 
%
\bibitem[Mr91]{Mather91} J. N. Mather 
{\em Action minimizing invariant measures for positive definite  
Lagrangian systems} Math. Z. {\bf 207}, 169-207 (1991).
%
\bibitem[Mr93]{Mather93} J. N. Mather 
{\em   Variational construction of connecting orbits}  Ann. Inst. Fourier (Grenoble)  43  (1993),  no. 5, 1349--1386.
%
\bibitem[Mr02]{Mather02}
J. N.  Mather {\em A property of compact, connected, laminated subsets of manifolds}  Ergodic Theory Dynam. Systems  22  (2002),  no. 5, 1507--1520.
%
\bibitem[Mr09]{Mather_snowbird}
J. N.  Mather {\em Order structure on action-minimizing orbits} Snowbird proceedings
%
\bibitem[MR95]{McShane-Rivin}
G. McShane, I. Rivin {\em A norm on homology of surfaces and counting simple geodesics}  Internat. Math. Res. Notices  1995,  no. 2, 61--69 
%
\bibitem[Mo86]{Moser} J. Moser 
{\em Minimal solutions of a variational problem on a torus} Ann. Inst. Henri Poincar\'e\ \textbf{3} (1986), 229--272 
%
\bibitem[O05]{Osvaldo} O. Osuna,  {\em Vertices of Mather's beta function}  
Ergodic Theory Dynam. Systems 25 (2005), no. 3, 949--955.
%
\bibitem[O09]{Osvaldo_09} O. Osuna {\em  The Aubry set for periodic Lagrangians on the circle }  Bol. Soc. Brasil. Mat. (N.S.)   40  (2009),  no. 2,  247-252.
%
\bibitem[S91]{Senn91} W. Senn {\em Strikte Konvexit\H{a}t f\H{u}r Variationsprobleme auf dem $n$-dimensionalen Torus}   Manuscripta Math.  71  (1991),  no. 1, 45--65.
%
\bibitem[S95]{Senn95} W. Senn {\em Differentiability properties of the minimal average action}  Calc. Var. Partial Differential Equations  3  (1995),  no. 3, 343--384.
%
\bibitem[S]{Sorrentino} A. Sorrentino {\em On the integrability of Tonelli Hamiltonians} preprint
%
\end{thebibliography} }
\end{document}